\documentclass[review,hidelinks,onefignum,onetabnum]{siamart250211}

\input{abbrev}

\usepackage{lipsum}
\usepackage{amsfonts}
\usepackage{graphicx}
\usepackage{epstopdf}
\usepackage{algorithmic}
\ifpdf
  \DeclareGraphicsExtensions{.eps,.pdf,.png,.jpg}
\else
  \DeclareGraphicsExtensions{.eps}
\fi

\newsiamremark{remark}{Remark}
\newsiamremark{hypothesis}{Hypothesis}
\crefname{hypothesis}{Hypothesis}{Hypotheses}
\newsiamthm{claim}{Claim}

\headers{Randomized Krylov methods for inverse problems}{Chung and Gazzola}

\title{Randomized Krylov methods for inverse problems\thanks{Submitted to the editors \today.
\funding{This work was partially supported by the National Science Foundation under grant DMS-2411197 (J. Chung). Any opinions, findings, conclusions or recommendations expressed in this material are those of the author(s) and do not necessarily reflect the views of the National Science Foundation.}
}}

\author{Julianne Chung\footnote{Department of Mathematics, Emory University, Atlanta, GA, 30322, USA}
\and Silvia Gazzola\footnote{Department of Mathematics, University of Pisa, 56126 Pisa, Italy}}
\usepackage{amsopn}

\usepackage{mathtools}

\usepackage{nicefrac}

\newcommand{\bfThetan}{\bfTheta^{(n)}}
\newcommand{\bfThetanp}{\bfTheta^{(m)}}

\newcommand{\tbfq}{\widetilde{\bfq}}

\newcommand{\tbfs}{\widetilde{\bfs}}

\newcommand{\reg}{\text{\scriptsize reg}}
\newcommand{\rec}{\text{\scriptsize rec}}
\newcommand{\true}{\text{\scriptsize true}}
\newcommand{\wtT}{\widetilde{\bfT}_{k+1}}

\newcommand{\Cmv}{C_{\text{\small{mv}}}}

\newcommand{\sk}{{\text{\small{sk}}}}

\ifpdf
\hypersetup{
  pdftitle={Randomized Krylov for inverse problems},
  pdfauthor={Julianne Chung and Silvia Gazzola}
}
\fi

\begin{document}
\nolinenumbers
\maketitle

\begin{abstract}
In this paper we develop randomized Krylov subspace methods for efficiently computing regularized solutions to large-scale linear inverse problems. Building on the recently developed randomized Gram-Schmidt process, where sketched inner products are used to estimate inner products of high-dimensional vectors, we propose a randomized Golub-Kahan approach that works for general rectangular matrices.  We describe new iterative solvers based on the randomized Golub-Kahan approach and show how they can be used for solving inverse problems with rectangular matrices, thus extending the capabilities of the recently proposed randomized GMRES method. 
We also consider hybrid projection methods that combine iterative projection methods, based on both the randomized Arnoldi and randomized Golub-Kahan factorizations, with Tikhonov regularization, where regularization parameters can be selected automatically during the iterative process. Numerical results from image deblurring and seismic tomography show the potential benefits of these approaches.
\end{abstract}



\section{Introduction}
Randomized approaches, in particular randomization for obtaining low-rank approximations or solving linear least-squares problems, have emerg\-ed as powerful tools in the numerical linear algebra community; see \cite{Martinsson_Tropp_2020} and references therein. Randomized Krylov methods are iterative methods that incorporate randomization in various ways: their main benefit is 
the potential of dramatically reducing 
computational costs with respect to their deterministic counterparts, while maintaining strong theoretical justifications and properties. The broad goal of this work is to develop randomized Krylov methods for large-scale linear inverse problems, which may be used as purely iterative regularization methods or in a hybrid fashion, i.e., combining iterative and variational regularization. 


Let us consider linear inverse problems of the form,
\begin{equation}
\label{eq:inverseproblem}
\bfb = \bfA \bfx_\true + \bfe
\end{equation}
where $\bfb\in \bbR^m$ collects observed measurements, $\bfA \in \bbR^{m\times n}$ represents the discretized forward model, $\bfx_\true \in \bbR^n$ is an unknown quantity of interest, 
and $\bfe\in \bbR^m$ contains noise or errors in the data.  Given $\bfb$ and $\bfA,$ the goal of the inverse problem is to approximate $\bfx_\true$. This is typically achieved by applying some kind of regularization to the original linear system \eqref{eq:inverseproblem}. Since, in general, $\bfA$ is large and unstructured, a common approach to recover $\bfx_\true$ is to apply an iterative (e.g., standard Krylov) solver to the least-squares problem
\begin{equation}
\label{eq:LS}
\min_{\bfx\in\bbR^n} \| \bfA \bfx - \bfb \|,
\end{equation}
where $\|\cdot\|$ denotes the vector 2-norm, sometimes also referred to as the $l_2$ norm. Early termination of many iterative solvers for \eqref{eq:LS} results in a regularized solution; see  \cite{Hansen2010}.  Another common approach to approximate $\bfx_\true$ is to solve a Tikhonov-regularized problem, i.e., compute
\begin{equation}
\label{eq:Tik}
\bfx_\reg =\argmin_{\bfx\in\bbR^n} \| \bfA \bfx - \bfb\|^2 + \lambda^2 \| \bfx \|^2=\argmin_{\bfx\in\bbR^n}\left\|
    \left[\begin{array}{cc}
    \bfA\\
    {\lambda}\bfI_n
    \end{array}\right]\bfx -
    \left[\begin{array}{cc}
    \bfb\\
    \bfzero
    \end{array}\right]\right\|^2,
\end{equation}
where $\lambda>0$ is the regularization parameter.  Many iterative (e.g., standard Krylov) methods can be used to solve minimization problems \cref{eq:Tik} efficiently, but they can get computationally expensive for large-scale problems that require many iterations.  Moreover, the regularization parameter $\lambda$ is typically not known a priori, and estimating it can be costly.  By combining Tikhonov regularization \eqref{eq:Tik} and iterative projection methods, hybrid projection methods provide a robust and computationally efficient approach for simultaneous solution approximation and regularization parameter selection; see  \cite{chung2024computational} for a survey paper with more details.

Randomized methods have been considered for solving linear least-squares problems like \cref{eq:LS,eq:Tik}. 
In this setting there are various ways to incorporate randomization, the most common ones being the so-called ``sketch-and-solve" and ``sketch-and-precondition" paradigms. In the first instance and considering an overdetermined matrix $\bfA$, sketching is initially applied to reduce the dimensionality of the matrix $\bfA$ and the data vector $\bfb$; a standard algorithm is then used to solve the approximate problem. Such approaches can provide low-rank matrix approximations, but they are not suited for inverse problems with a slowly decaying spectrum. In the context of inverse problems, a randomized approach was proposed in \cite{landman2025randomized}, where the inner product free LSLU projection method was combined with ``sketch-and-solve'' on the projected problem. In the second instance, a sketching matrix is used to obtain a preconditioner for the considered iterative methods: for example, this underlies the LSNR \cite{lsnr} and the Blendenpik \cite{Avron2010Blendenpik} algorithms.  There is a vast body of literature (often within the statistical and data science remit), which explores randomization to solve regularized linear least-squares problems, formulated as damped least-squares problems; see, for instance, \cite{7836598, RidgeSketch, pmlr-v80-chowdhury18a,meier2022randomizedalgorithmstikhonovregularization}. Although these sketching-based strategies have gained widespread interest as tools for low-rank matrix approximation and for solving linear systems, they are numerically unstable for ill-conditioned, and even moderately ill-conditioned, problems; see \cite{Meier2024Stable}. Moreover, their utility for state-of-the-art inverse problems remains largely under-explored.
Apart from the efficient solution of linear least-squares problems, other previous works on randomized Krylov methods focus on block Krylov methods for spectral computation and matrix approximation (e.g., for approximating extreme singular values and singular vectors); see, e.g., \cite{musco2015randomized, tropp2022randomized}. More recently, randomized iterative solvers have been developed, where sketching techniques are used to reduce the computational cost of (re)orthogonalization \cite{BalabanovRGMRES2022}, and have been used to speed up approximating matrix functions \cite{cortinovis2024speeding}. Specifically, a randomized Gram-Schmidt process was described in \cite{BalabanovRGMRES2022}, which forms the basis for the randomized Arnoldi (rArnoldi) algorithm and the randomized GMRES (rGMRES) method. These approaches are restricted to square matrices.

In this work, we develop randomized Krylov methods for large-scale inverse problems where randomization is incorporated via sketched inner products.  The main contributions of this work are the following. 
\begin{itemize}
    \item  Extending the approach of \cite{BalabanovRGMRES2022}, we develop a randomized Golub-Kahan (rGK) algorithm that can be applied to a rectangular matrix.  This approach requires two different sketching matrices and generates an upper Hessenberg and an upper triangular matrix. We provide cost and storage comparisons to (standard) Golub-Kahan bidiagonalization.
    \item Then, building on the proposed rGK approach, we introduce various randomized iterative projection methods, where the main difference among the methods is the constraint imposed on the residuals or normal equations residuals. Specifically, we describe a randomized LSQR (rLSQR) method, a randomized CGLS (rCGLS) method, and a randomized LSMR (rLSMR) method.
    \item For all of the considered randomized Krylov methods (rGMRES for square systems, and rLSQR, rCGLS, and rLSMR for both square and rectangular systems), we describe their hybrid versions, where we estimate the regularization parameter during the iterative process by applying well-established regularization parameter selection methods.
\end{itemize}

An outline for this paper is as follows.  In \Cref{sec:randomized}, we review the rArnoldi and rGMRES algorithms for square matrices, and we propose the rGK algorithm for rectangular matrices.  We describe various randomized iterative methods for the linear least-squares problem \eqref{eq:LS} based on the rGK factorization, and draw connections to their standard counterparts. Hybrid projection methods based on these randomized projection methods and Tikhonov regularization are described in \Cref{sec:hybridrandomized}.  Numerical results are provided in \Cref{sec:numerics}, and conclusions and future work are described in \Cref{sec:conclusions}.

\section{Randomized Krylov methods}
\label{sec:randomized}
Krylov subspace methods have been widely studied and used for solving large-scale inverse problems, where each iteration requires a matrix-vector product with $\bfA$ and, possibly, $\bfA\t$.  Various works have reduced the cost of matrix-vector products, e.g., by exploiting structure or sparsity in $\bfA$ or utilizing highly parallel distributed architectures; see \cite{NaPaPe04, van2015astra} just to cite a few.  However, a remaining challenge is the computational cost of orthonormalization. In particular, computing inner products can be expensive if the vectors are very high-dimensional or the number of iterations is high.  Moreover, inner products can present a communication bottleneck if vectors are distributed across multiple processors. Modifications and variants such as inner product free iterative methods \cite{brown2024hlslu,brown2024hcmrh} have been considered, but often lead to quasi-minimum residual methods and can require expensive pivoting strategies. 

Randomized approaches have been proposed within Krylov  methods to reduce the cost of orthonormalization; see, e.g., \cite{BalabanovRGMRES2022, NakatsukasaSGMRES2024}. Most of these methods rely on the following random sketching technique to estimate inner products.  Let $\calV$ be a subspace of $\bbR^n$ of low dimension, and let $\bfThetan \in \bbR^{\ell_n \times n}$, with $\ell_n \ll n$ be a sketching matrix. The $l_2$ inner product between vectors in subspaces of $\bbR^n$ is estimated by
\begin{equation}
\langle \cdot, \cdot \rangle \approx \langle \bfThetan \cdot, \bfThetan \cdot \rangle.
\end{equation}
More precisely, we say that $\bfThetan$ is an $\eps$ embedding for the subspace $\calV$ if
\begin{equation}\label{def:epsemb}
\left\vert\langle\bfx,\bfy\rangle-\langle\bfThetan\bfx,\bfThetan\bfy\rangle\right\vert\leq\eps\|\bfx\|\|\bfy\|,\quad\mbox{for all $\bfx,\bfy\in\calV$}.     
\end{equation}
The above relation implies that
\begin{equation}\label{prop:epsemb}
    (1+\eps)^{-1}\|\bfThetan\bfx\|^2\leq \|\bfx\|^2\leq (1-\eps)^{-1}\|\bfThetan\bfx\|^2.
\end{equation}
Adopting a randomized approach allows for savings in that inner products of high dimensional vectors are replaced by inner products of reduced vectors obtained by sketching. This is of course efficient as far as the random sketching matrix can be cheaply applied to a vector. For instance, $\bfThetan$ can be chosen as a Gaussian matrix, although applying it to a vector costs $O(\ell_n n)$ flops; when $\bfThetan$ is chosen as the more sophisticated Subsampled Randomized Hadamard Transform (SRHT), applying it to a vector costs as little as $2n\log_2(\ell_n+1)$ flops, depending on the implementation and the computational architecture; see \cite{SRHT}. We refer to \cite{Woodruff} for a survey of other possible sketching matrices. Since, in the remaining part of the paper, the subspace $\calV$ is a Krylov subspace generated as the iterations of the considered solver unfold, we need to consider sketching matrices that do not require a-priori knowledge of $\calV$ except for its dimension. More precisely, we consider $(\eps,\delta,K+1)$-oblivious embeddings such that \eqref{def:epsemb} is satisfied with probability at least $1-\delta$ for any fixed $(K+1)$-dimensional subspace of $\bbR^n$. Both Gaussian and SRHT matrices have this property for big enough $\ell_n$.

To establish notation and as motivation, we begin in \Cref{sub:rGMRES} with a review of the randomized Gram-Schmidt process that was described in \cite{BalabanovRGMRES2022}, and the resulting rArnoldi and rGMRES algorithms. Then, in \Cref{sub:rLSQR}, we propose a new randomized Golub-Kahan (rGK) approach for rectangular matrices, where two different sketching matrices are used to efficiently project the problem using two different Krylov subspaces. Various randomized Krylov solvers based on the rGK approach are described in \Cref{sub:randsolvers}. 

\subsection{rArnoldi and rGMRES}
\label{sub:rGMRES}

Recall that standard GMRES and the underlying Arnoldi algorithm can be used to solve high-dimensional inverse problems of the form \cref{eq:inverseproblem}, where $\bfA \in \bbR^{n \times n}$ and $\bfb \in \bbR^n$; see, for instance, \cite{JeHa07}. By taking the first basis vector $\bar{\bfq}_1$ as the normalized right-hand-side vector $\bfb$, $k\leq K$ iterations of the Arnoldi algorithm generate a matrix with orthonormal columns $\bar \bfQ_{k+1}=[\bar\bfq_1,\dots,\bar\bfq_{k+1}] \in \bbR^{n \times (k+1)}$ and an upper Hessenberg matrix 
$\bar \bfH_k \in \bbR^{(k+1) \times k}$ that satisfy the Arnoldi identity
\begin{equation}
\label{eq:Arnoldi} \bfA \bar \bfQ_{k} = \bar \bfQ_{k+1} \bar \bfH_k.
\end{equation}
The columns $\bar{\bfq}_i$, $i=1,\dots,k+1$, of $\bar \bfQ_{k+1}$ span the Krylov subspace $\calK_{k+1}(\bfA,\bfb) = {\rm span}\{ \bfb, \bfA \bfb, \ldots, \bfA^k \bfb\}$, i.e., 
$\calR(\bar\bfQ_{k+1}) = \calK_{k+1}(\bfA,\bfb)$. Here and in the following, $\calR(\cdot)$ denotes the range or column space of a matrix. At the $k$th iteration of the GMRES method, the solution to \cref{eq:LS} can be approximated by
\begin{equation}
\label{eq:GMRES_projectedproblem}
\bar\bfx_k = \bar\bfQ_k \bar\bfz_k\,, \quad \mbox{where} \quad \bar\bfz_k = \argmin_\bfz \|\bar\bfH_k \bfz - \bar \beta \bfe_1 \|\quad\mbox{with}\quad \bar \beta=\|\bfb\|\,.
\end{equation}
Here and in the following, $\bfe_i$ denotes the $i$th canonical basis vector, whose size should be clear from the context. 
The Arnoldi algorithm can be regarded as a specific instance of the Gram-Schmidt process, which is used to construct an orthonormal basis for $\calK_k(\bfA,\bfb) $ in an iterative process. 
Note that the orthonormalization process can be costly for very large dimensional vectors.  

One remedy to mitigate the computational cost of the Gram-Schmidt process is to use the randomized Gram-Schmidt process that exploits sketching, as described in \cite{BalabanovRGMRES2022}. Given a matrix $\bfW\in\bbR^{n\times K}$ and an $\ell_n \times n$ sketching matrix $\bfThetan$, with $\ell_n \ll n$, two matrices $\bfQ \in \bbR^{n \times K}$ and $\bfR \in \bbR^{K \times K}$ are generated, 
where the columns of $\bfQ$ are orthogonal with respect to the sketched inner product $\langle \bfThetan \cdot, \bfThetan \cdot \rangle$ and $\bfR$ is an upper triangular matrix.  More specifically, at the $k$th step, $2\leq k\leq K$, we compute the strictly upper triangular entries $\bfr_{k-1}$ of $\bfR$ as
\begin{equation}\label{rGS_step}
\bfr_{k-1}=(\underbrace{\bfThetan\bfQ_{k-1}}_{=:\bfS_{k-1}})\t(\bfThetan\bfw_k),
\end{equation}
where $\bfw_k$ is the $k$th column of $\bfW$ and $\bfQ_{k-1} \in \bbR^{n \times (k-1)}$ contains the first $(k-1)$ columns of $\bfQ$. This is followed by an update for $\bfq_k$ of the form
\[
\bfq_k=\nicefrac{(\bfw_k-\bfQ_{k-1}\bfr_{k-1})}{\|\bfThetan(\bfw_k-\bfQ_{k-1}\bfr_{k-1})\|}.
\]
In \eqref{rGS_step} we use the notation $\bfS_{k-1} := \bfThetan \bfQ_{k-1} \in \bbR^{\ell_n \times {(k-1)}}$ 
to emphasize the reduction in the size of the vectors that are orthonormalized. Note that, in \cite{BalabanovRGMRES2022}, \eqref{rGS_step} is replaced by the more general least-squares problem,
\[
\bfr_{k-1}=\argmin_{\bfy\in\bbR^{k-1}} \| \bfThetan \bfQ_{k-1} \bfy - \bfThetan \bfw_k \|.
\]
Indeed, the vector $\bfr_{k-1}$ in \eqref{rGS_step} is a specific solution expression corresponding to the classical Gram-Schmidt implementation. Although more sophisticated ways of computing $\bfr_{k-1}$ (such as stable solvers based on Householder transformations or Richardson iterations for the associated normal equations) allow for greater stability of the randomized Gram-Schmidt process, we did not observe substantial differences when testing on a variety of inverse problems, so here we use the simpler version \eqref{rGS_step}. 

Now, the randomized Gram-Schmidt process can be used to compute an approximation of $\bfx_\true$ in \eqref{eq:inverseproblem} in the subspaces $\calK_k(\bfA,\bfb)$, $k\leq K$. Let $\bfThetan$ be an $(\varepsilon,\delta,K+1)$ oblivious subspace embedding (implying that it is an $\varepsilon$-embedding for any $\calK_k(\bfA, \bfb)$, $k\leq K+1$). The 
randomized Arnoldi (rArnoldi) approach is obtained by applying the randomized Gram-Schmidt algorithm to the Krylov subspace $\calK_k(\bfA,\bfb)$. More precisely, after $k$ iterations of rArnoldi, we have a basis $\{\bfq_1,\dots,\bfq_{k+1}\}$ of $\calK_{k+1}(\bfA, \bfb)$ 
with columns $\bfq_i$, $i=1,\dots,k+1$, that are orthonormal with respect to the sketched inner product $\langle \bfThetan \cdot, \bfThetan \cdot \rangle$, and an upper Hessenberg matrix $\bfH_k \in \bbR^{(k+1) \times k}$. Let 
$\bfQ_{k}=[\bfq_1,\dots,\bfq_{k}]\in\bbR^{n\times k}$ and $\bfQ_{k+1}=[\bfQ_k,\bfq_{k+1}]\in\bbR^{n\times (k+1)}$, then we get an Arnoldi identity analogous to \cref{eq:Arnoldi}, i.e.,
\begin{equation}
\label{eq:rArnoldi} \bfA  \bfQ_{k} = \bfQ_{k+1} \bfH_k
\end{equation}
in infinite precision arithmetic. Building on the rArnoldi algorithm, the rGMRES solution at the $k$th iteration is given by
\begin{equation}
\label{eq:rGMRES_projectedproblem}
\bfx_k = \bfQ_k \bfz_k\,,\quad\mbox{where}\quad\mbox\bfz_k = \argmin_\bfz \|\bfH_k \bfz - \beta \bfe_1 \|\quad\mbox{and}\quad \beta=\|\bfThetan \bfb\|\,.
\end{equation}
In infinite precision, the $k$th rGMRES approximation $\bfx_k$ \cref{eq:rGMRES_projectedproblem} minimizes the sketched (semi)norm of the residual over $\calK_k(\bfA,\bfb)$, since
\begin{equation}
\label{eq:GMRES_sketchedprojected}
\begin{split}
\|\bfH_k\bfz_k-\beta\bfe_1\|&=\min_\bfz \|\bfH_k \bfz - \beta \bfe_1 \| =
\min_\bfz \|\bfThetan \bfQ_{k+1} (\bfH_k \bfz - \beta \bfe_1 ) \|\\ 
&= \min_\bfz \|\bfThetan (\bfA \bfQ_k \bfz - \bfb) \|=\|\bfThetan (\bfA \bfx_k-\bfb)\|\,.
\end{split}
\end{equation}
This optimality property of rGMRES can be also derived from more general facts about projection methods 
(as described, for instance, in \cite[Chapter 5]{saad2003iterative}): indeed, 
the $k$th rGMRES approximation is such that 
\[
\bfx_k\in\calK_k(\bfA,\bfb)\quad\mbox{and}\quad \bfr_k\perp_{\bfThetan}\bfA\calK_k(\bfA,\bfb),
\]
where, by $\perp_{\bfThetan}$, we mean orthogonal with respect to the sketched inner product $\langle\bfThetan\cdot,\bfThetan\cdot\rangle$. Note that the above properties generalize the ones used by standard GMRES, which would hold with $\bfThetan=\bfI_n$. Using the fact that $\bfThetan$ is an $(\eps,\delta,K+1)$-oblivious embedding, we can prove the following bound for the rGMRES residual $\bfb-\bfA\bfx_k$ norm with respect to the GMRES residual $\bfb-\bfA\bar\bfx_k$ norm at the $k$th iteration
\[
\|\bfb-\bfA\bfx_k\|^2\leq \frac{(1+\eps)}{(1-\eps)}\|\bfb-\bfA\bar\bfx_k\|^2,
\]
which comes directly from the inequalities
\begin{eqnarray*}
(1-\eps)\|\bfb-\bfA\bfx_k\|^2\leq
\|\bfThetan(\bfb-\bfA\bfx_k)\|^2\leq
\|\bfThetan(\bfb-\bfA\bar\bfx_k)\|^2\leq
(1+\eps)\|\bfb-\bfA\bar\bfx_k\|^2,
\end{eqnarray*}
where we have used \eqref{prop:epsemb} for the first and third inequality, and the rGMRES optimality property \eqref{eq:GMRES_sketchedprojected} for the second inequality.

The randomized GMRES algorithm is provided in \Cref{alg:rGMRES}. 
Note that, with respect to the version proposed in \cite{BalabanovRGMRES2022} and in addition to the simplified orthogonalization procedure in step \ref{alg:orthogonal}, step \ref{alg:skecthq} replaces the computation of $\tbfs_{k+1}$ with $\tbfs_{k+1}= \bfThetan \tbfq_{k+1}$. It can be easily shown by strong induction that the two expressions are equivalent. If applying sketching to an $n$-dimensional vector costs more than $\ell_n(k+1)$ flops, then our implementation is more efficient than the one in \cite{BalabanovRGMRES2022}. 


\begin{algorithm}
    \begin{algorithmic}[1]
        \REQUIRE{$\bfA\in\bbR^{n \times n}$, $\bfb\in\bbR^n$, $\bfThetan\in\bbR^{\ell_n \times n}$ ($\ell_n \ll n$)}.
    \STATE Initialize $\tbfq_1 = \bfb$.
    \STATE Sketch $\tbfs_1 = \bfThetan \tbfq_1$.
    \STATE Compute the sketched (semi)norm $\beta=\|\tbfs_1\|$.
    \STATE Scale vectors $\bfs_1=\tbfs_1/\beta$, $\bfq_1=\tbfq_1/\beta$.
    \FOR{$k= 1,2, \dots K$}
    \STATE {\textbf{Randomized Arnoldi}}
    \STATE Get new vector $\tbfq_{k+1} = \bfA \bfq_k$.
    \STATE\label{alg:LSsolve} Sketch $\tbfs_{k+1} = \bfThetan \tbfq_{k+1}$.
\STATE\label{alg:orthogonal} Compute $[\bfR]_{(1:k,k+1)} = \bfS_{k}\t\tbfs_{k+1}$.
\STATE Compute projection: $\tbfq_{k+1} = \tbfq_{k+1} - \bfQ_{k} [\bfR]_{(1:k,k+1)}$.
\STATE\label{alg:skecthq} Compute the sketched projection 
$\widetilde{\bfs}_{k+1}=\widetilde{\bfs}_{k+1} - \bfS_{k}[\bfR]_{(1:k,k+1)}$.
\STATE Compute the sketched (semi)norm $r_{k+1,k+1} = \| \tbfs_{k+1} \|$. 
\STATE Scale vectors: $\bfs_{k+1} = \tbfs_{k+1} / r_{k+1,k+1}$, 
$\bfq_{k+1} = \tbfq_{k+1} / r_{k+1,k+1}$.
\STATE {\textbf{Solve the projected problem}}
$$
\bfz_k = \argmin_\bfz \|\bfH_k \bfz - \beta \bfe_1 \| ,\quad
\mbox{where}\quad \bfH_k=[\bfR]_{(1:k+1,2:k+1)}.
$$
\STATE Compute approximate solution $\bfx_k = \bfQ_k \bfz_k$.
    \ENDFOR
    \end{algorithmic}
    \caption{rGMRES method}
    \label{alg:rGMRES}
\end{algorithm}

\subsection{Randomized Golub-Kahan Approach}
\label{sub:rLSQR}
The aim of this section is to introduce a randomized Golub-Kahan (rGK)
algorithm that works for general rectangular matrices $\bfA$. First recall that, given $\bfA\in\bbR^{m\times n}$ and $\bfb\in\bbR^{m}$, $k$ iterations of the standard Golub-Kahan bidiagonalization (GKB) algorithm generate two matrices $\bar \bfU_{k+1} \in \bbR^{m \times (k+1)}$ and $\bar \bfV_{k} \in \bbR^{n \times k}$ with numerically orthonormal columns such that 
\begin{equation}
    \label{eq:2Krylovsubspaces}
\calR(\bar\bfV_k)=\calK_k(\bfA\t\bfA,\bfA\t\bfb)\subset \bbR^n\quad\mbox{and}\quad\calR(\bar\bfU_{k+1})=\calK_{k+1}(\bfA\bfA\t,\bfb)\subset \bbR^{m},
\end{equation}
i.e., the columns of $\bar\bfV_k$ and $\bar\bfU_{k+1}$ span the Krylov subspaces detailed above. The matrices $\bar\bfV_k$ and $\bar\bfU_{k+1}$ are linked by the relations
\begin{equation}\label{eq: sketchip}
    \bfA \bar \bfV_k = \bar \bfU_{k+1} \bar \bfB_{k} \quad\mbox{and}\quad \bfA\t \bar \bfU_{k+1} = \bar \bfV_{k+1} (\bar \bfB_{k+1,k})\t, 
\end{equation}
where $\bar\bfB_{k} \in \bbR^{(k+1) \times k}$ is lower bidiagonal and $\bar\bfB_{k+1,k}$ is obtained by removing the last column of $\bar\bfB_{k+1}$.

Similar to the randomized approaches described in \Cref{sub:rGMRES}, the rGK method is obtained by replacing standard inner products appearing in GKB with sketched inner products.  Since, in general, $m\neq n$, we need to consider two different random sketchings for estimating the inner products in the two Krylov subspaces in  \eqref{eq:2Krylovsubspaces}.  Assume that we have identified a maximum number of iterations $K$ such that $k\leq K$. Let $\bfTheta^{(n)} \in \bbR^{\ell_n \times n}$, with $\ell_n \ll n$, be an $(\eps,\delta,K+1)$-oblivious sketching matrix for vectors in $\bbR^n$, and let  $\bfTheta^{(m)} \in \bbR^{\ell_m \times m}$, with $\ell_m \ll m$, be an $(\eps,\delta,K+1)$-oblivious sketching matrix for vectors in $\bbR^m$.
That is, the $l_2$ inner product between vectors in subspaces of dimensions up to $K$ in $\bbR^n$ and $\bbR^{m}$ are estimated by
\begin{equation}
\label{eq:sketchrGK}
\langle \cdot, \cdot \rangle \approx \langle \bfTheta^{(n)} \cdot, \bfTheta^{(n)} \cdot \rangle\quad\mbox{and}\quad
\langle \cdot, \cdot \rangle \approx \langle \bfTheta^{(m)} \cdot, \bfTheta^{(m)} \cdot \rangle,
\end{equation}
respectively. If $m=n$, 
one can simply take $\bfTheta^{(n)}=\bfTheta^{(m)}$. 

A first observation is that, while the traditional GKB algorithm 
enjoys short recurrences, its randomized version will require long (i.e., full Gram-Schmidt-like) recurrences. To understand why this is the case, in a simpler setting, one can consider generating an orthonormal basis for the subspace $\calK_k(\widehat{\bfA}, \bfb)$, where $\widehat{\bfA}\in\bbR^{n\times n}$ is symmetric and $\bfb\in\bbR^n$. It is well known that, in this setting, the Arnoldi algorithm reduces to the symmetric Lanczos algorithm, whereby the leading $k\times k$ minor of the upper Hessenberg matrix $\bar{\bfH}_k$ is a symmetric tridiagonal matrix. One way of motivating why this happens is that, directly from \eqref{eq:Arnoldi} (with $\bfA=\widehat{\bfA}$), 
\begin{equation}\label{eq:symmTri}
\bar\bfq_j\t\widehat{\bfA}\bar\bfq_i=h_{j,i} \mbox{ ($=0$ if $j>i+1$)}\quad\mbox{and}\quad
\bar\bfq_i\t\widehat{\bfA}\bar\bfq_j=h_{i,j} \mbox{ ($=0$ if $i>j+1$)}\,.
\end{equation}
Since $\widehat{\bfA}$ is symmetric, $\bar\bfq_j\t\widehat{\bfA}\bar\bfq_i=\bar\bfq_i\t\widehat{\bfA}\bar\bfq_j$ and, therefore, $h_{i,j}=h_{j,i}=0$ unless \linebreak[4]$i=j-1,j,j+1$. In a randomized setting, when generating a $\bfTheta^{(n)}$-orthonormal basis for $\calK_k(\widehat{\bfA},\bfb)$ we still get a relation of the kind \eqref{eq:Arnoldi}, but relations of the kind \eqref{eq:symmTri} do not hold anymore because, due to the lack of orthogonality of the columns of $\bfQ_{k+1}=[\bfq_1,\dots,\bfq_{k+1}]\in\bbR^{n\times (k+1)}$,
\[
\bfq_i\t\widehat{\bfA}\bfq_j\neq\bfq_j\t\widehat{\bfA}\bfq_i=\bfq_j\t\bfQ_{k+1}{\bfH}_{k}\bfe_i\neq \bfzero\quad\text{for $i,j=1,\dots,k$,}
\]
in general. Indeed, in this setting, one has a relation of the kind
\[
\bfQ_{k+1}\t((\bfTheta^{(n)})\t)\bfTheta^{(n)}\widehat{\bfA}\bfQ_k={\bfH}_k\,,
\]
so that the leading $k\times k$ minor of ${\bfH}_k$ is not symmetric unless the matrices $(\bfThetan)\t\bfThetan$ and $\widehat{\bfA}$ commute (which is very unlikely). 


After $k\leq K$ iterations of the rGK algorithm, we have the following two partial matrix factorizations
\begin{equation}\label{eq: radGK}
\bfA\bfV_k = \bfU_{k+1}\bfM_k\,,\quad\mbox{and}\quad
\bfA\t\bfU_{k+1} = \bfV_{k+1}\bfT_{k+1}\,,
\end{equation}
where 
\begin{itemize}
    \item $\bfV_k=[\bfv_1,\dots,\bfv_k]\in\bbR^{n\times k}$ and $\bfV_{k+1}=[\bfV_k,\bfv_{k+1}]\in\bbR^{n\times (k+1)}$ have $\bfThetan$-orthogonal columns
    \item $\bfU_{k+1}=[\bfu_1,\dots,\bfu_{k+1}]\in\bbR^{m\times (k+1)}$ has $\bfThetanp$-orthogonal columns
    \item $\bfM_k\in\bbR^{(k+1)\times k}$ is an upper Hessenberg matrix
    \item $\bfT_{k+1}\in\bbR^{(k+1)\times (k+1)}$ is an upper triangular matrix
\end{itemize}
and the sketched vectors are the (orthonormal) columns of the matrices
\begin{equation}
\bfP_k = \begin{bmatrix}\bfp_1 & \ldots \bfp_k \end{bmatrix}  \in \bbR^{\ell_n \times k}
\quad \mbox{and} \quad \bfS_k = \begin{bmatrix}\bfs_1 & \ldots \bfs_k \end{bmatrix} \in \bbR^{\ell_m \times k} .
\end{equation}
Algorithm \ref{alg: randGK} summarizes the main steps involved in the new rGK algorithm. 
\begin{algorithm}
    \begin{algorithmic}[1]
        \REQUIRE{$\bfA\in\bbR^{m \times n}$, $\bfb\in\bbR^{m}$,  $\bfThetan\in\bbR^{\ell_n \times n}$ ($\ell_n \ll n$), $\bfThetanp\in\bbR^{\ell_m \times m}$ ($\ell_m \ll m$).}
    \STATE Initialize $\widetilde{\bfu}_1 = \bfb$.
    \STATE Sketch $\widetilde{\bfs}_1 = \bfThetanp \widetilde{\bfu}_1 $.
    \STATE Compute the sketched (semi)norm 
    $\beta=\| \widetilde{\bfs}_1 \|$.
    \STATE Scale vectors $\bfs_1=\widetilde{\bfs}_1/\beta$, $\bfu_1=\widetilde{\bfu}_1/\beta$.
    \STATE Compute $\widetilde{\bfv}_1 = \bfA\t\bfu_1$.
    \STATE Sketch $\widetilde{\bfp}_1 = \bfThetan\widetilde{\bfv}_1$.
    \STATE Compute the sketched (semi)norm $t_{1,1}=\| \widetilde{\bfp}_1 \|$.
    \STATE Scale vectors $\bfp_1=\widetilde{\bfp}_1 / t_{1,1}$, $\bfv_1=\widetilde{\bfv}_1/t_{1,1}$. 
    \FOR{$k= 1, 2, \dots K$}
    \STATE {\textbf{Randomized Gram-Schmidt for constructing $\bfU$}}
    \STATE Get new vector $\widetilde{\bfu}_{k+1} = \bfA\bfv_k$.
    \STATE Sketch $\widetilde{\bfs}_{k+1} = \bfThetanp \widetilde{\bfu}_{k+1}$.
     \STATE\label{alg:LSsolveGK1} Compute 
$[\bfM]_{(1:k,k)} = \bfS_{k}\t\widetilde{\bfs}_{k+1}$.
\STATE Compute the projection: $\widetilde{\bfu}_{k+1} = \widetilde{\bfu}_{k+1} - \bfU_{k}[\bfM]_{(1:k,k)}$.
\STATE\label{alg:supdateGK} Compute the sketched projection: $\widetilde{\bfs}_{k+1} = \widetilde{\bfs}_{k+1} - \bfS_{k}[\bfM]_{(1:k,k)}$. 
\STATE Compute the sketched (semi)norm $m_{k+1,k} = \| \widetilde{\bfs}_{k+1} \|$.
\STATE Scale vectors: $\bfs_{k+1}= \widetilde{\bfs}_{k+1}/ m_{k+1,k}$,
$\bfu_{k+1}  = \widetilde{\bfu}_{k+1} / m_{k+1,k}$.
\smallskip
\STATE {\textbf{Randomized Gram-Schmidt for constructing $\bfV$}}
    \STATE Get new vector $\widetilde{\bfv}_{k+1} = \bfA\t\bfu_{k+1}$.
    \STATE\label{alg:pupdate} Sketch $\widetilde{\bfp}_{k+1} = \bfThetan\widetilde{\bfv}_{k+1}$.
    \STATE\label{alg:LSsolveGK2} Compute 
$[\bfT]_{(1:k,k+1)} = \bfP_{k}\t\widetilde{\bfp}_{k+1}$.
\STATE Compute the projection: $\widetilde{\bfv}_{k+1} = \widetilde{\bfv}_{k+1} - \bfV_{k} [\bfT]_{(1:k,k+1)}$.
\STATE\label{alg:pupdateGK} Compute the sketched projection 
$\widetilde{\bfp}_{k+1} = \widetilde{\bfp}_{k+1}-\bfP_k[\bfT]_{(1:k,k+1)}$.
\STATE Compute the sketched (semi)norm $t_{k+1,k+1} = \| \widetilde{\bfp}_{k+1} \|$.
\STATE Scale vectors: $\bfp_{k+1} = \widetilde{\bfp}_{k+1} / t_{k+1,k+1}$, 
$\bfv_{k+1} = \widetilde{\bfv}_{k+1} / t_{k+1,k+1}$.
    \ENDFOR
    \end{algorithmic}
    \caption{rGK method}\label{alg: randGK}
\end{algorithm}
Note that, similarly to what was done for rArnoldi in Algorithm \ref{alg:rGMRES}, steps \ref{alg:LSsolveGK1} and \ref{alg:LSsolveGK2} of Algorithm \ref{alg: randGK} can be replaced by the more general tasks of approximating the solution to the least-squares problems,
\[
[\bfM]_{(1:k,k)} = \argmin_{\bfy} \| \bfS_{k} \bfy -\widetilde{\bfs}_{k+1}\|
\quad\mbox{and}\quad [\bfT]_{(1:k,k+1)} = \argmin_{\bfy} \| \bfP_{k} \bfy -\widetilde{\bfp}_{k+1} \|,
\]
respectively. Also, depending on the cost of performing sketching, one may replace
steps \ref{alg:supdateGK} and \ref{alg:pupdateGK} of Algorithm \ref{alg: randGK} by the sketched vectors
\[
\widetilde{\bfs}_{k+1} = \bfThetanp\widetilde{\bfu}_{k+1}\quad\mbox{and}\quad\widetilde{\bfp}_{k+1} = \bfThetan\widetilde{\bfv}_{k+1},
\]
respectively.

We now present an estimate of the cost, in terms of flops, of computing the rGK factorization \eqref{eq: radGK}, and compare it to the classical GKB. Let us denote by $\Cmv$ the cost of matrix-vector products with $\bfA$ and let us assume this equals the cost of matrix-vector products with $\bfA\t$; let us denote by $C_{\sk_n}$ and $C_{\sk_m}$ the cost of performing matrix-vector products with the sketching matrices $\bfTheta^{(n)}$ and $\bfTheta^{(m)}$, respectively. Then the cost of ($K$ iterations of) Algorithm \ref{alg: randGK} is estimated to be
\begin{equation}\label{rGKcost}
\begin{split}
(2K+1)\Cmv &+ (K+1)(m + n +C_{\sk_n}+C_{\sk_m}+ 3(\ell_n+\ell_m) - 2)\\
&+K(K+1)(m + n + 2\ell_m + 2\ell_n -1)
\end{split}
\end{equation}
flops. In contrast, the cost of ($K$ iterations of) the (standard) 
GKB 
algorithm is
\[
(2K + 1)\Cmv + (K+1)(5(m + n) -2) -2(m+n)
\]
flops. If 
GKB is implemented with 
reorthogonalization (ro-GKB), then its computational cost is
\[
(2K+1)\Cmv + (K+1)(5(m + n) -2) -2(m+n) + \frac{K(K-1)}{2}(3(n + m)- 2)
\]
flops. 
A few remarks are in order. All the above estimates have the term \linebreak[4]$(2K+1)\Cmv$ in common, accounting for the matrix-vector multiplications with $\bfA$ and $\bfA\t$ performed to expand the underlying Krylov subspaces (up to dimension $(K+1)$). The terms of the kind $O(K(m+n))$ account for rescaling and sums of (full-dimensional) vectors of length $m$ and $n$. In rGK and 
ro-GKB, the terms of the kind $O(K^2(m+n))$ account for the computations of the bases of the underlying Krylov subspaces (up to dimension $(K+1)$), and they appear with different scaling factors; in particular, the count for 
ro-GKB also includes scalar products computations. In rGK, the terms of the kind $O(K(\ell_n + \ell_m))$, and $O(K^2(\ell_n+\ell_m))$, account for rescaling and sums, and inner products, of sketched vectors, respectively.   
More precisely, when it comes to the cost of performing inner products, 
rGK requires
\[
(K+1)(K+2)(\ell_m+\ell_n-1)\,\, \mbox{flops;}
\]
ro-GKB requires
\[
K(K+1)(n + m -1)\,\, \mbox{flops;}
\]
standard GKB requires the computations of two norms at each iteration, totaling
\[
2K(n+m-1)\,\, \mbox{flops.}
\]
Therefore, in terms of inner product computations and recalling that $\ell_m\ll m$ and $\ell_n\ll n$, rGK is more convenient than 
ro-GKB when $C_{\sk_m}$ and $C_{\sk_n}$ are negligible. 
%

\begin{figure}
    \centering
    \begin{tabular}{cc}
    \small{\textbf{(a)}} & \small{\textbf{(b)}}\vspace{-0.1cm}\\
    \includegraphics[width=5cm]{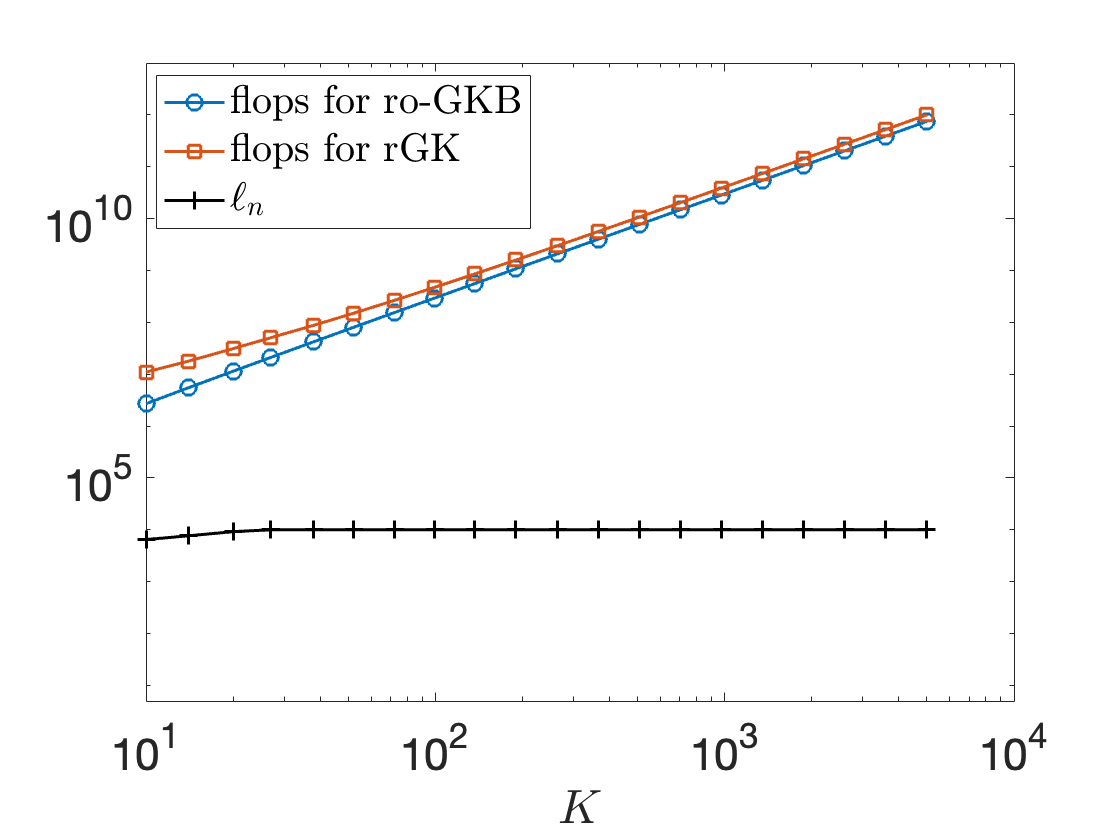}& 
    \includegraphics[width=5cm]{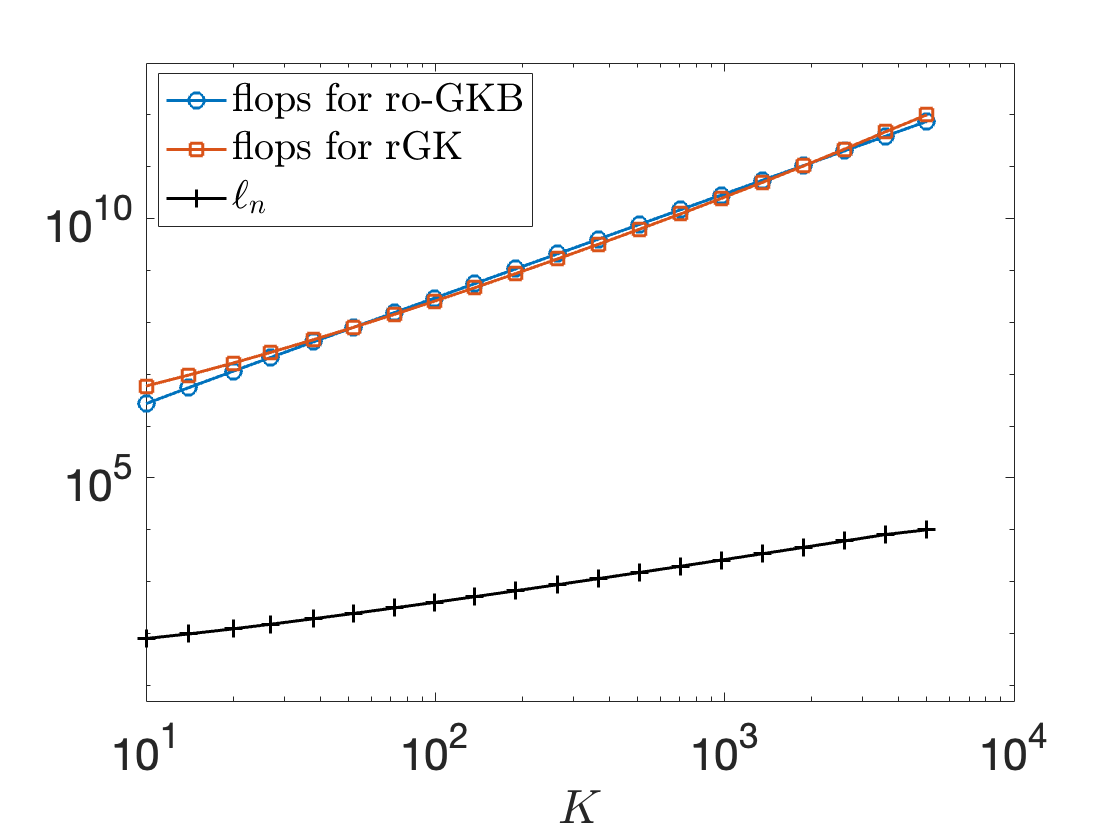}\\
    \small{\textbf{(c)}} & \small{\textbf{(d)}}\vspace{-0.1cm}\\
    \includegraphics[width=5cm]{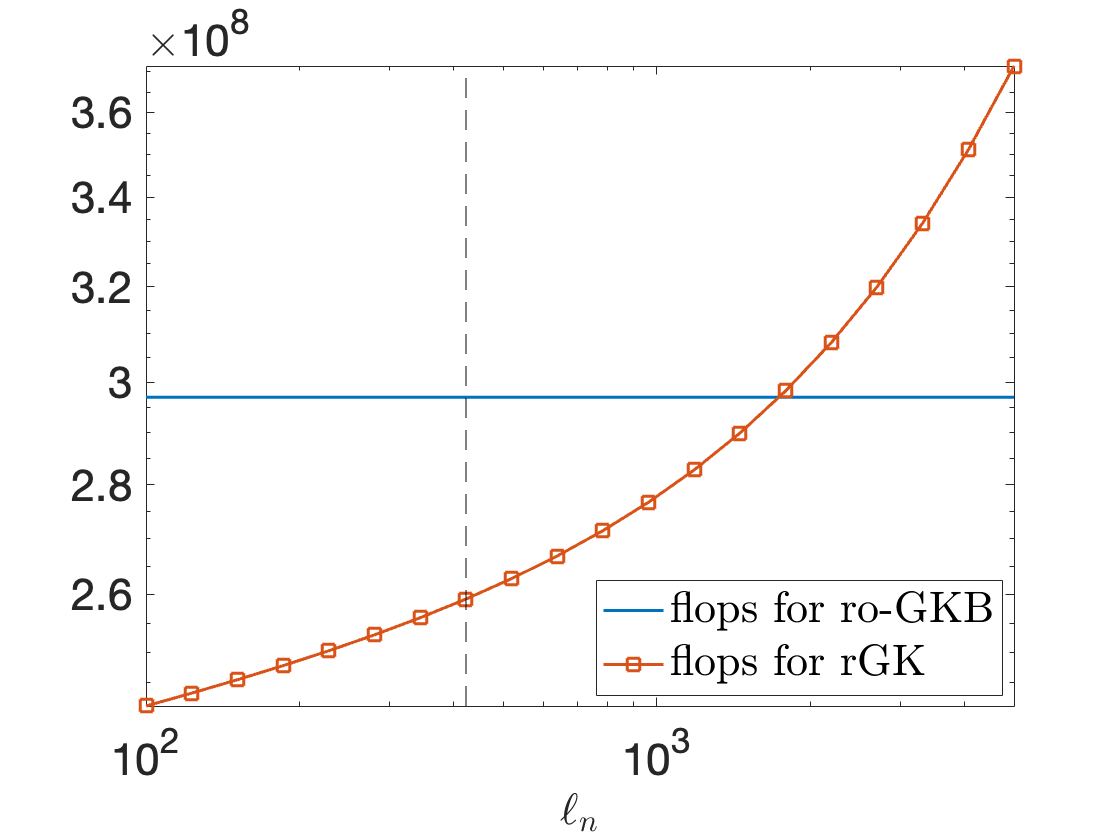} & 
    \includegraphics[width=5cm]{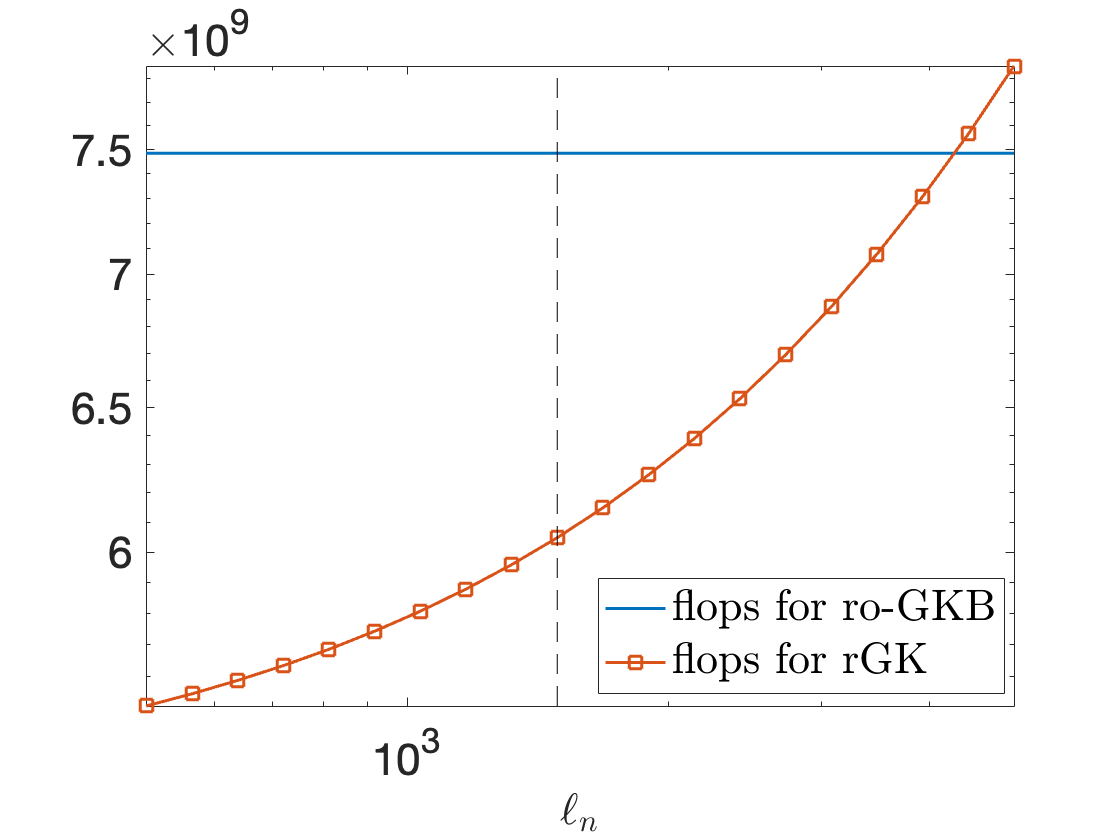}
    \end{tabular}
    \caption{Computational cost comparisons of ro-GKB and rGK, for a matrix of size $m=n=10^4$. Flops and $\ell_n$ against 20 approximately logarithmically equispaced values of $K$ between 10 and 5000, when $\ell_n$ is chosen according to \eqref{ellnest} (frame \textbf{(a)}) and according to \eqref{ellnesth} (frame \textbf{(b)}). Flops against 20 approximately logarithmically equispaced values of $\ell_n$ between $K$ and 5000, with $K=100$ (frame \textbf{(c)}) and $K=500$ (frame \textbf{(d)}). In the bottom frames, a vertical dashed line corresponds to the value of $\ell_n$ selected according to \eqref{ellnesth}.
    }
    \label{fig:flops_comp}
\end{figure}

To support the above reasoning about computational costs, we consider \linebreak[4]$m=n=10^4$ and we take $\bfTheta^{(m)}=\bfTheta^{(n)}$ as a SRHT (also mentioned at the beginning of \Cref{sec:randomized}). 
It is well-known (see, e.g., \cite{BalabanovRGMRES2022}) that SRHT is an $(\eps, \delta, K+1)$ oblivious $l_2$-subspace embedding if
\begin{equation}\label{ellnest}
\ell_n\geq 2(\eps^2-\nicefrac{\eps^3}{3})^{-1}\left(\sqrt{K+1} + \sqrt{8\log(\nicefrac{6n}{\delta})}\right)^2\log(\nicefrac{3(K+1)}{\delta}).
\end{equation}
In Figure \ref{fig:flops_comp} we compare the estimated flops required to run rGK and ro-GKB when $\ell_n$ is selected according to \eqref{ellnest}, with $\eps=\delta=0.5$, as well as with the more heuristic choice
\begin{equation}\label{ellnesth}
\ell_n\geq{\nicefrac{2(K+1)\log(n)}{\log(K+1)}},
\end{equation}
and when 20 approximately logarithmically equispaced values of $K$ between 10 and 5000 are considered. We also run comparisons when $K$ is fixed and 20 approximately logarithmically equispaced values of $\ell_n$ between $K$ and 5000 are selected. We can see that, when $\ell_n$ is kept relatively low, there is a clear computational gain in performing rGK rather than ro-GKB. 
We conclude this subsection by mentioning that, when regularizing inverse problems, both in a variational and Bayesian setting, it is known that ro-GKB may perform better than GKB, as the former may deliver more accurate estimations of the regularized solution itself or of the so-called hyperparameters defining it (like the Tikhonov regularization parameter $\lambda$ in \eqref{eq:Tik}), often requiring less iterations than standard GKB; see  \cite{baglama1996iterative,calvetti2000ribbon,saibaba2020efficient} for more details. Therefore, in this setting, one should adopt ro-GKB, and the new rGK may provide clear computational advantages over it.

\subsection{Randomized rGK-based solvers for \cref{eq:LS}} 
\label{sub:randsolvers}
Building upon the rGK algorithm, a number of solvers can be introduced to approximate a solution to problem \cref{eq:LS}. The solvers vary based on the conditions imposed on the residual.

We start by deriving a new \emph{randomized LSQR} (rLSQR) method that, at the $k$th iteration, defines the $k$th approximate solution as
\begin{equation}\label{eq: projrandLSQR}
\bfx_k=\bfV_k\bfz_k\in\calK_k(\bfA\t\bfA,\bfA\t\bfb),\quad \mbox{where}\quad 
\bfz_k=\argmin_{\bfz\in\bbR^k}\|\bfM_k \bfz - \beta\bfe_1\|\,,    
\end{equation}
and where $\bfV_k$ and $\bfM_k$ are the outputs of $k$ iterations of \Cref{alg: randGK}, where also $\beta(=\|\bfThetanp\bfb\|)$ is defined. In exact arithmetic, it can be easily seen that $\bfx_k$ minimizes the norm of the sketched residual over the subspace $\calK_k(\bfA\t\bfA, \bfA\t\bfb)$. Indeed,
\begin{align}
    \|\bfM_k\bfz_k - \beta \bfe_1\|=\min_{\bfz\in\bbR^k}\|\bfM_k\bfz - \beta \bfe_1\| &=
    \min_{\bfz\in\bbR^k}\|\bfThetanp\bfU_{k+1}(\bfM_k\bfz - \beta \bfe_1)\|\nonumber\\
&=\min_{\bfz\in\bbR^k}\|\bfThetanp(\bfA\bfV_k\bfz - \bfb)\|\,,\label{rLSQRoptprop}
\end{align}
where, in the second equality, we have exploited the orthonormality of the columns of $\bfThetanp\bfU_{k+1}$, and, in the third equality, we have used the first of the factorizations in \eqref{eq: radGK}.  
These derivations are analogous to the ones performed to define rGMRES starting from the rArnoldi process (see \eqref{eq:GMRES_sketchedprojected}); the optimality property for these two methods is analogous, too. 
In the general framework for iterative projection methods for \eqref{eq:LS} (see again \cite[Chapter 5]{saad2003iterative}), rLSQR can be regarded as a projection method that, at the $k$th iteration, computes
\begin{equation}\label{rLSQR_conds}
\bfx_k\in\calK_k(\bfA\t\bfA,\bfA\t\bfb)\quad\mbox{such that}\quad \bfr_k\perp_{\bfTheta^{(m)}}\bfA\calK_k(\bfA\t\bfA,\bfA\t\bfb),
\end{equation}
where $\perp_{\bfThetanp}$ denotes orthogonality with respect to the second sketched inner product in \eqref{eq:sketchrGK}. 
This generalizes LSQR, for which the constraints on the residual are given by $\bfr_k\perp\bfA\calK_k(\bfA\t\bfA,\bfA\t\bfb)$. According to the general theory of projection methods, 
\eqref{rLSQR_conds} implies that the $\bfTheta^{(m)}$-seminorm of $\bfr_k$ is minimized at each iteration, which agrees with 
\eqref{rLSQRoptprop}. Eventually, one could also derive rLSQR starting from the normal equations associated to \eqref{rLSQRoptprop}, i.e., by equaling the gradient of the objective function in \eqref{rLSQRoptprop} to zero, i.e., 
\[
\bfA\t(\bfTheta^{(m)})\t\bfTheta^{(m)}\bfA\bfx = \bfA\t(\bfTheta^{(m)})\t \bfTheta^{(m)} \bfb\,.
\]
The $k$th rLSQR approximate solution is such that 
\begin{equation}\label{rLSQR_conds_var}
\bfx_k\in\calK_k(\bfA\t\bfA,\bfA\t\bfb)\quad\mbox{and}\quad
\bfA\t(\bfTheta^{(m)})\t\bfTheta^{(m)}\bfr_k\perp\calK_k(\bfA\t\bfA,\bfA\t\bfb).
\end{equation}
Note that the standard $l_2$ inner product is used to impose the orthogonality condition above. 
Exploiting the rGK factorizations \eqref{eq: radGK} we recover the projected problem
\[
\bfM_k\t\bfM_k\bfz_k = \beta\bfM_k\t\bfe_1\,,\quad
\mbox{and take}\quad \bfx_k=\bfV_k\bfz_k,
\]
which is equivalent to \eqref{eq: projrandLSQR}. 
For $\bfThetanp=\bfI_m$, i.e., in a deterministic setting, conditions \eqref{rLSQR_conds} and \eqref{rLSQR_conds_var} coincide with the ones imposed in a deterministic setting to derive the mathematically equivalent LSQR and CGLS methods, respectively.

We next derive a rGK-based solver starting from the conditions
\begin{equation}\label{rCGLS_conds}
\bfx_k\in\calK_k(\bfA\t\bfA,\bfA\t\bfb)\quad\mbox{and}\quad
\bfA\t\bfr_k\perp_{\bfThetan}\calK_k(\bfA\t\bfA,\bfA\t\bfb)\,,
\end{equation}
where $\perp_{\bfThetan}$ denotes orthogonality with respect to the first sketched inner product in \eqref{eq:sketchrGK}. Exploiting the rGK factorizations \eqref{eq: radGK}, imposing the above conditions amounts to computing
\begin{equation}\label{eq: projrandCGLS}
\bfx_k=\bfV_k\bfz_k,\quad\mbox{where $\bfz_k\in\bbR^k$ is such that $\wtT\bfM_k\bfz_k=\beta t_{1,1}\bfe_1$}\,,
\end{equation}
and where $\bfV_k$ and $\bfM_k$ are obtained running $k$ iterations of Algorithm \ref{alg: randGK}, \linebreak[4]$\wtT=[\bfI_k, \bfzero]\bfT_{k+1}\in\bbR^{k\times (k+1)}$ is $\bfT_{k+1}$ without its last row, $t_{1,1}$ is the (1,1) entry of $\bfT_{k+1}$, and $\beta=\|\bfThetanp\bfb\|$. 
For $\bfThetan=\bfI_n$ and $\bfThetanp=\bfI_m$, i.e., in a deterministic setting, conditions \eqref{rCGLS_conds} would be equivalent to \eqref{rLSQR_conds_var}, which would underlie the derivation of CGLS (mathematically equivalent to LSQR). However, because of randomness, \eqref{rCGLS_conds} leads to the solver \eqref{eq: projrandCGLS} that is not mathematically equivalent to rLSQR. Because of the analogy with the deterministic setting, we regard $\bfx_k$ in  
\eqref{eq: projrandCGLS} as the $k$th approximate solution computed by the new 
\emph{randomized CGLS} (rCGLS) method. Note, however, that rCGLS does not bear other  resemblances with standard CGLS, 
including the lack of short-recurrences for the solution updates. Equivalently, one could derive rCGLS starting from the equations
\[
(\bfThetan)\t\bfThetan\bfA\t\bfA\bfx=(\bfThetan)\t\bfThetan\bfA\t\bfb\,,
\]
and imposing similar conditions to the ones in \eqref{rCGLS_conds}, but with orthogonality of the residual in the standard $l_2$ inner product. 
Straightforward derivations involving the rGK partial factorizations \eqref{eq: radGK} show that one obtains the same projected problem as in \eqref{eq: projrandCGLS}. The lack of optimality properties of this solver is evident (also) from the equation above, as the coefficient matrix $(\bfThetan)\t\bfThetan\bfA\t\bfA$ is not even symmetric (unless $(\bfThetan)\t(\bfThetan)$ and $\bfA\t\bfA$ commute, which is very unlikely). Indeed, in this setting, none of the well-known optimality results for general iterative projection methods apply; see \cite[Chapter 5]{saad2003iterative}.

We finally derive a new \emph{randomized LSMR} (rLSMR) method  that, at the $k$th iteration, defines the $k$th approximate solution of \eqref{eq:LS} by imposing the conditions 
\begin{equation}\label{rLSMR_conds}
\bfx_k\in\calK_k(\bfA\t\bfA,\bfA\t\bfb)\quad\mbox{and}\quad
\bfA\t\bfr_k\perp_{\bfThetan}\bfA\t\bfA\calK_k(\bfA\t\bfA,\bfA\t\bfb)\,.
\end{equation}
Running $k\leq K$ iterations of Algorithm \ref{alg: randGK} and taking $\bfx=\bfV_k\bfz$, with $\bfz\in\bbR^k$, the second condition above and the rGK factorizations \eqref{eq: radGK} lead to the linear system of equations
\[
(\bfThetan\bfA\t\bfA\bfV_k)\t(\bfThetan\bfA\t\bfU_{k+1}(\beta\bfe_1-\bfM_k\bfz))=\bfzero\,,
\]
which can be further simplified to compute the $k$th rLSMR solution as
\[
\bfx_k=\bfV_k\bfz_k,\quad\mbox{where}\quad 
\bfM_k\t\bfT_{k+1}\t(\beta t_{1,1}\bfe_1-\bfT_{k+1}\bfM_k\bfz_k)=\bfzero,
\]
or, equivalently,
\begin{equation}\label{rLSMR_projpb}
\bfx_k=\bfV_k\bfz_
k,\quad\mbox{where}\quad 
\bfz_k=\argmin_{\bfz\in\bbR^k}\|\bfT_{k+1}\bfM_k\bfz-\beta t_{1,1}\bfe_1\|;
\end{equation}
also here, $\beta=\|\bfThetanp\bfb\|$ and $t_{1,1}$ is the (1,1) entry of $\bfT_{k+1}$. In exact arithmetic, it can be easily seen that rLSMR minimizes the norm of the sketched normal equations residual over the subspace $\calK_k(\bfA\t\bfA,\bfA\t\bfb)$. Indeed,
\begin{align}
    \|\bfT_{k+1}\bfM_k\bfz_k-\beta t_{1,1}\bfe_1\|&=\min_{\bfz\in\bbR^k}\|\bfT_{k+1}\bfM_k\bfz-\beta t_{1,1}\bfe_1\| \nonumber\\
&=\min_{\bfz\in\bbR^k}\|\bfThetan\bfV_{k+1}\bfT_{k+1}(\bfM_k\bfz-\beta \bfe_1)\|\nonumber\\
&=\min_{\bfz\in\bbR^k}\|\bfThetan\bfA\t\bfU_{k+1}(\bfM_{k}\bfz-\beta \bfe_1)\|\nonumber\\
&=\min_{\bfx\in\calK_k(\bfA\t\bfA,\bfA\t\bfb)}\|\bfThetan\bfA\t(\bfA\bfx-\bfb)\|.\label{rLSMRoptprop}
\end{align}
We note that both conditions \eqref{rLSMR_conds} and the optimality property \eqref{rLSMRoptprop} generalize to the randomized setting the standard LSMR ones, which would hold with $\bfThetan=\bfI_n$.


Although the new rGK-based solvers introduced in this section are 
natural generalizations of their classical counterparts, their implementations may significantly differ. Most notably, as a consequence of the lack of short recurrences for the updates of the basis vectors $\bfv_i$, $\bfu_i$, $i=1,\dots,k+1$ (and lack of bidiagonal structure in the matrices $\bfM_k$ and $\bfT_{k+1}$), the entire basis should be stored as the iterations proceed, so that the storage cost increases linearly with the number of iterations. Moreover, the sketched basis vectors collected in the matrices $\bfP_{k+1}\in\bbR^{\ell_n\times (k+1)}$ and $\bfS_{k+1}\in\bbR^{\ell_m\times (k+1)}$ should be stored (although, if $\ell_n\ll n$ and $\ell_m\ll m$, their overhead is minimal). These are the trade-offs of having an $l_2$ inner product free and randomized method, but in the inverse problems setting, they are mitigated by the following facts. First of all, as mentioned at the end of \Cref{sub:rLSQR}, in many situations it may be more appropriate to run ro-GKB rather than standard GKB, demanding the storage of all $\bfv_i$ and $\bfu_i$ vectors anyway. Second, when performing (standard) GKB-based hybrid projection methods (described for the randomized case in the next section), in order to allow adaptive regularization parameter choice strategies at each iteration, one has to store the solution basis vectors $\bfv_i$.  We also mention that similar storage requirements occur in flexible Golub-Kahan methods to handle Tikhonov-regularized problems more general than \eqref{eq:Tik}; see \cite{chung2019flexible}. Such memory demands may be mitigated, also in the randomized setting, by appropriately using some recycling strategies; see also \cite{jiang2021hybrid}.  

\paragraph{SVD approximation} A common way to assess the regularization capabilities of iterative projection methods 
is to study the quality of the approximation of the dominating singular values of $\bfA$ obtained by taking the singular values of the corresponding projected matrices; see, for instance, \cite{GazzolaSilvia2016Iotd, gazzola2015survey, Hansen2010}. Although a theoretical investigation of such properties for the new rGK-based solvers is outside the scope of this work, we provide an illustration to empirically show the different performance of such solvers.

\begin{figure}
    \centering
    \begin{tabular}{ccc}
\small{$\sigma(\bfM_k)$} & 
\small{$\sqrt{\sigma(\widetilde{\bfT}_{k+1}\bfM_k)}$} & 
\small{$\sqrt{\sigma({\bfT_{k+1}}\bfM_k)}$}\\
\includegraphics[width=0.3\linewidth]{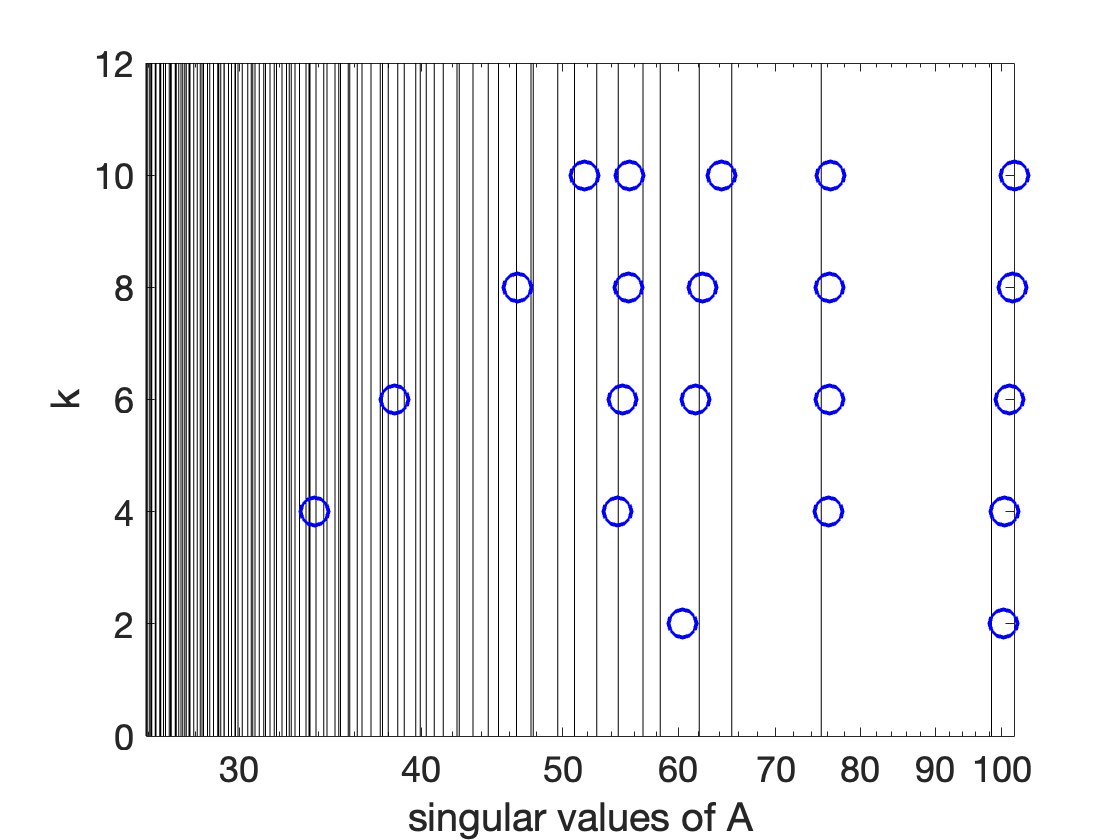}    & \includegraphics[width=0.3\linewidth]{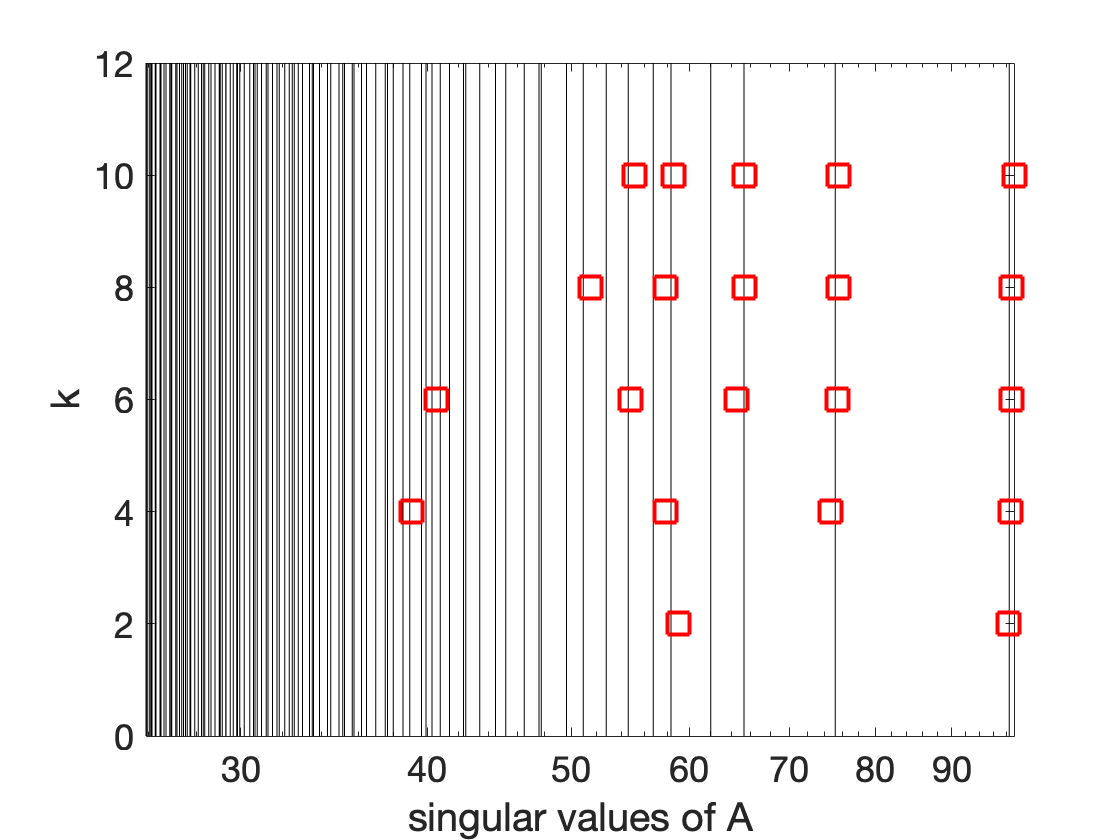}
         & \includegraphics[width=0.3\linewidth]{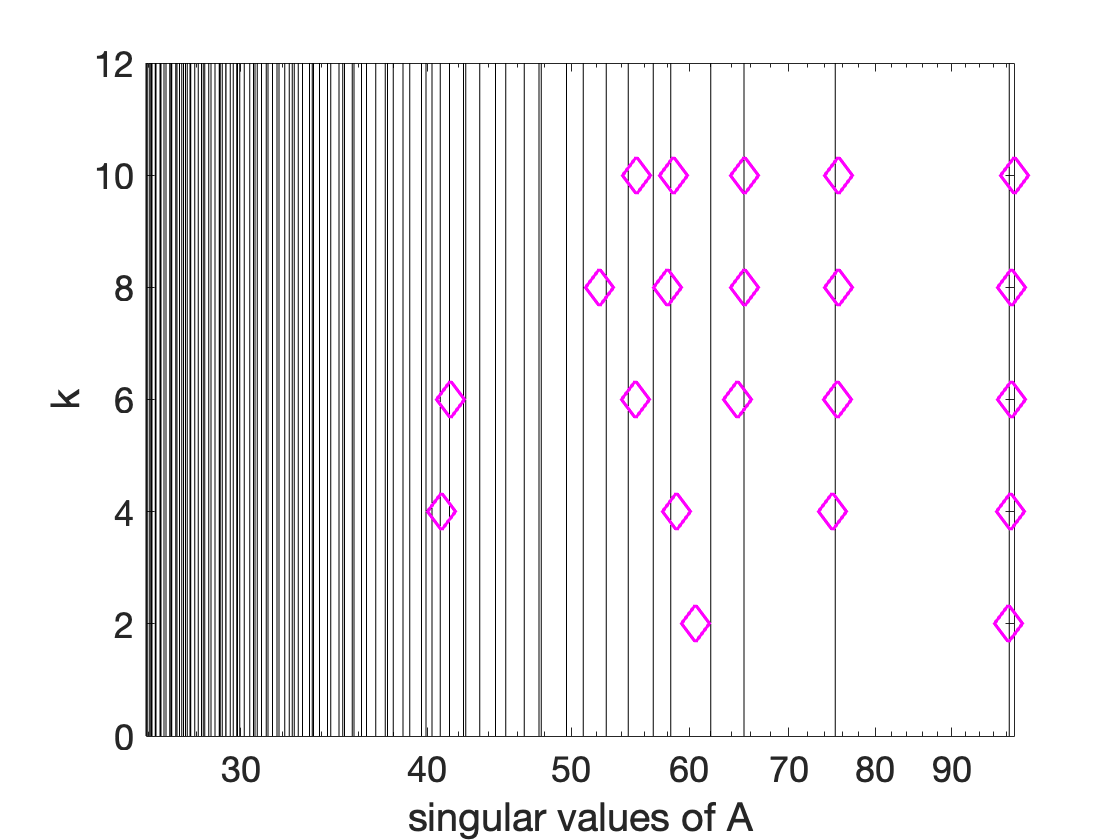}
    \end{tabular}
    \caption{Approximation of the dominating singular values of $\bfA$ (solid vertical lines) via the singular values of the projected matrices associated to rLSQR (leftmost frame), rCGLS (middle frame) and rLSMR (rightmost frame), at iterations $k=2,4,6,8,10$.}
    \label{fig:Ex3_SVDapprox}
\end{figure}

%
%
We generate a sparse matrix $\bfA$ of size $8192\times 4096$ modelling a seismic tomography problem (see \Cref{sec:numerics}, Example 2, for a description of a larger instance of this test problem). 
%
In all the frames in \cref{fig:Ex3_SVDapprox} we plot the dominating singular values of $\bfA$ with solid vertical lines. Then, for $k=2,4,6,8,10$, $k$ different markers highlight the singular values of the $k$th projected matrices associated to different  solvers (we show only the ones falling within the dominating singular values of $\bfA$). As $k$ increases, we observe that the rCGLS and rLSMR projected matrices provide better approximations than rLSQR to the larger singular values of $\bfA$. This may affect the behavior of these solvers when employed as iterative regularization methods, and when employed in a hybrid fashion (see the next section), especially when it comes to choosing the regularization parameter for the latter. 

\section{Hybrid randomized Krylov methods} 
\label{sec:hybridrandomized}
In this section, we describe how the solvers in \Cref{sub:rGMRES} and \Cref{sub:randsolvers} can be applied to approximate the solution of the standard Tikhonov regularized problem \cref{eq:Tik}, with possible adaptive choice of the regularization parameter $\lambda\geq 0$, in the framework of hybrid projection methods, as outlined in  \cite{chung2024computational}. For most of this section, we will naturally consider the first formulation in \eqref{eq:Tik}.

For standard Krylov subspace methods, one starting point to introduce hybrid projection methods, which also justifies them from a theoretical point of view, is the fact that first projecting the original problem and then regularizing the projected problem is equivalent to first regularizing the original problem and then projecting the regularized problem; see \cite{Hansen2010}. In the following we show that, when considering the randomized Krylov subspace methods introduced in the previous subsections, this is only true when considering specific instances of using the sketched seminorms $\|\bfThetanp\cdot\|$ and $\|\bfThetan\cdot\|$. We assume, for now, that $\lambda$ is fixed. 

We begin by considering projection methods based on the rArnoldi factorization \eqref{eq:rArnoldi}, and we show that regularizing the projected problem is not equivalent to projecting the regularized problem \eqref{eq:Tik}.  That is, consider projecting the LS problem,
\[
\min_{\bfx\in\bbR^n}\|\bfA\bfx-\bfb\|^2
\]
onto the Krylov subspace $\calK_k(\bfA,\bfb)$, $k\leq K$,
which amounts to solving
\[
\min_{\bfz\in\bbR^k}\|\bfV_{k+1}({\bfH}_k\bfz-\beta\bfe_1)\|^2,\quad\mbox{with}\quad  \beta=\|\bfThetan\bfb\|.
\]
We denote the project-then-regularize solution $\bfx_{k,\lambda}:=\bfV_k\bfz_{k,\lambda}\in\calK_k(\bfA,\bfb)$, where
\begin{equation}\label{eq: FPTR}
\bfz_{k,\lambda}  = \argmin_{\bfz\in\bbR^k}\|\bfV_{k+1}({\bfH}_k\bfz-\beta\bfe_1)\|^2+\lambda^2\|\bfz\|^2\,.
\end{equation}
On the other hand, if we start with the Tikhonov problem \eqref{eq:Tik} and project it using rArnoldi (i.e., taking 
$\bfx \in \calR(\bfV_k)$ in
\eqref{eq:Tik}), we get the regularize-then-project solution \linebreak[4]$\bfx_{\lambda,k}:=\bfV_k\bfz_{\lambda,k}\in\calK_k(\bfA,\bfb)$, where
\begin{equation}\label{eq: FRTP}
\bfz_{\lambda,k} = \argmin_{\bfz\in\bbR^k}\|\bfV_{k+1}({\bfH}_k\bfz-\beta\bfe_1)\|^2+\lambda^2\|\bfV_k\bfz\|^2,\quad\mbox{with}\quad \beta=\|\bfThetan\bfb\|.
\end{equation} 
Although both provide approximations to the regularized solution $\bfx_\reg$ in \eqref{eq:Tik}, it is clear that 
$\bfx_{\lambda,k}\neq\bfx_{k,\lambda}$, as problems \eqref{eq: FPTR} and \eqref{eq: FRTP} are not equivalent and cannot be further simplified (because $\bfV_{k+1}$ does not have orthonormal columns). However, since $\bfV_{k}$ has $\bfThetan$-orthonormal columns, replacing the regularization term in \eqref{eq:Tik} and \eqref{eq: FRTP} with the $\bfThetan$-seminorm results in the same projected problem,
    \[
    \min_{\bfz\in\bbR^{k}}\|\bfV_{k+1}(\bfH_{k}\bfz-\beta\bfe_1)\|^2+\lambda^2\|\bfz\|^2.
    \]
A further simplified projected problem can be obtained by replacing also the fit-to-data norm in \eqref{eq:Tik}, \eqref{eq: FPTR}, and \eqref{eq: FRTP} with the $\bfThetan$- seminorm. That is, the $k$th hybrid rGMRES iterate is given by $\bfx_k=\bfV_k\bfz_k$, where
\begin{equation}\label{eq: randreg}
\bfz_k = \argmin_{\bfz\in\bbR^{k}}\|\bfH_{k}\bfz-\beta\bfe_1\|^2+\lambda^2\|\bfz\|^2,
\end{equation}
and it holds that $\bfx_k$ is a minimizer of
    \begin{equation}
    \label{eq:projected}
    \min_{\bfx\in\calK_{k}(\bfA,\bfb)}\|\bfThetan(\bfA\bfx-\bfb)\|^2+\lambda^2\|\bfThetan\bfx\|^2\,,
    \end{equation}
which is equivalent to damped LS problem,
    \begin{equation}
    \label{eq:dampedprojected}
    \min_{\bfx\in\calK_{k}(\bfA,\bfb)}\left\|\left[
    \begin{array}{cc}
    \bfThetan & \\
    & \bfThetan 
    \end{array}\right]\left(
    \left[\begin{array}{cc}
    \bfA\\
    {\lambda}\bfI_n
    \end{array}\right]\bfx -
    \left[\begin{array}{cc}
    \bfb\\
    \bfzero
    \end{array}\right]\right)\right\|^2.
    \end{equation}
Moving to the rGK-based solvers, we can consider hybrid variants, where all the above rArnoldi derivations (including the lack of equivalence between the project-then-regularize and regularize-then-project approaches) generalize to this case. 
After running $k\leq K$ rGK iterations, we define the $k$th hybrid rLSQR iterate by $\bfx_k=\bfV_k\bfz_k$, where
\begin{equation}\label{eq:hybrid_rLSQR1}
\bfz_k=\argmin_{\bfz\in\bbR^k}\|\bfM_k\bfz-\beta\bfe_1\|^2+\lambda^2\|\bfz\|^2,
\end{equation}
and we have that 
\begin{equation}\label{eq:hybrid_rLSQR2}
\bfx_k=\argmin_{\bfx\in\calK_k(\bfA\t\bfA,\bfA\t\bfb)}\|\bfThetanp(\bfA\bfx-\bfb)\|^2+\lambda^2\|\bfThetan\bfx\|^2,
\end{equation}
or, equivalently,
\begin{equation}\label{eq:rLSQR_Tik}
    \bfx_k=\argmin_{\bfx\in\calK_{k}(\bfA\t\bfA,\bfA\t\bfb)}\left\|\left[
    \begin{array}{cc}
    \bfThetanp & \\
    & \bfThetan 
    \end{array}\right]\left(
    \left[\begin{array}{cc}
    \bfA\\
    {\lambda}\bfI_n
    \end{array}\right]\bfx -
    \left[\begin{array}{cc}
    \bfb\\
    \bfzero
    \end{array}\right]\right)\right\|^2.
\end{equation}
Alternatively, in the rGK setting, one could consider applying rLSQR directly to the second Tikhonov problem formulation in \eqref{eq:Tik}, i.e., Tikhonov formulated as a damped least-squares problem. At the $k$th rLSQR iteration, we get that 
\begin{equation}\label{eq:rLSQR_dampTik}
\bfx_k=\argmin_{\bfx\in\calK_k((\bfA\t\bfA+\lambda^2\bfI_n),\bfA\t\bfb)}\left\|
\bfTheta^{(m+n)}
   \left(
    \left[\begin{array}{cc}
    \bfA\\
    {\lambda}\bfI_n
    \end{array}\right]\bfx -
    \left[\begin{array}{cc}
    \bfb\\
    \bfzero
    \end{array}\right]\right)\right\|^2
\end{equation}
is an approximation to $\bfx_\reg$, where $\bfTheta^{(m+n)}\in\bbR^{\ell_{m+n}\times (m+n)}$ is a $(\eps, \delta, K+1)$-oblivious sketching matrix for vectors in $\bbR^{m+n}$. Note that problem \eqref{eq:rLSQR_Tik} can be regarded as a special case of \eqref{eq:rLSQR_dampTik}, as $\calK_k((\bfA\t\bfA+\lambda^2\bfI_n),\bfA\t\bfb)=\calK_k(\bfA\t\bfA,\bfA\t\bfb)$. However the advantage of the hybrid method defined in \eqref{eq:rLSQR_Tik} over the sketched Tikhonov formulation \eqref{eq:rLSQR_dampTik} is that the former naturally allows to set the regularization parameter $\lambda$ adaptively; this is not so straightforward when considering the latter, as the matrices $\bfM_k$ and $\bfT_{k+1}$ appearing in the rGK decomposition \eqref{eq: radGK} depend on $\lambda$; see also \cite{FrMa99}. 
    To the best of our knowledge, most of the other randomized methods for Tikhonov regularization consider a sketch-and-solve approach and are naturally formulated as \eqref{eq:rLSQR_dampTik} without the constraint on the approximate solution (i.e., taking $\bfx\in\bbR^n$), which is not necessarily computed through an iterative method. None of these solvers allow for easy tuning of the regularization parameter; see, e.g., \cite{RidgeSketch,7836598}.

One could also consider a hybrid version of rCGLS and rLSMR that incorporates Tikhonov regularization for approximating $\bfx_\reg$. Starting with rCGLS, at the $k$th iteration, $k\leq K$, condition \eqref{rCGLS_conds} can be applied where we replace $\bfA$ with the matrix $[\bfA\t, \lambda\bfI_n]\t$ and $\bfr_k$ with the vector $[(\bfb-\bfA\bfx_k)\t,\lambda\bfx_k\t]\t$, both naturally appearing in the second formulation of the Tikhonov problem \eqref{eq:Tik}. Then, using the properties of the matrices in the rGK factorization \eqref{eq: radGK}, we get the $k$th rCGLS approximation applied to the leftmost problem in \eqref{eq:Tik} as
\begin{equation}\label{eq:hybrid_rCGLS1}
\bfx_k=\bfV_k\bfz_k\in\calK_k(\bfA\t\bfA,\bfA\t\bfb)=\calK_k((\bfA\t\bfA+\lambda^2\bfI_n),\bfA\t\bfb),
\end{equation}
where $\bfz_k$ is such that
\begin{equation}\label{eq:hybrid_rCGLS2}
(\widetilde{\bfT}_{k+1}\bfM_k + \lambda^2\bfI_n)\bfz_k=\beta t_{1,1}\bfe_1.
\end{equation}
The above equation can be regarded as a regularization of the projected system appearing in \eqref{eq: projrandCGLS}, and we refer to the solver so obtained as hybrid rCGLS. Finally, one can regularize the rLSMR projected problem \eqref{rLSMR_projpb} by taking $\bfx_k=\bfV_k\bfz_k\in\calK_k(\bfA\t\bfA,\bfA\t\bfb)$, where
\begin{equation}\label{eq:hybrid_rLSMR1}
\bfz_k=\argmin_{\bfz\in\bbR^k}\|\bfT_{k+1}\bfM_k\bfz-\beta t_{1,1}\bfe_1\|^2+\lambda^2\|\bfz\|^2\,.
\end{equation}
For the hybrid rLSMR method so defined it holds that
\begin{equation}\label{eq:hybrid_rLSMR2}
\bfx_k=\argmin_{\bfx\in\calK_k(\bfA\t\bfA,\bfA\t\bfb)}\|\bfThetan\bfA\t(\bfA\bfx-\bfb)\|^2+\lambda^2\|\bfThetan\bfx\|^2,
\end{equation}
which can be regarded as a Tikhonov-regularized version of problem \eqref{rLSMRoptprop}. Contrary to hybrid rLSQR and hybrid rCGLS, hybrid rLSMR is not linked to the original Tikhonov formulation because the fit-to-data term is defined with respect to the normal equations residual. Indeed, even for (standard) LSMR, it is well-known that project-then-regularize is not equivalent to regularize-then-project; see \cite{chung2015hybrid} for more details. When rLSMR is concerned, we just briefly mention that one may obtain a different regularize-then-project rLSMR approximation of $\bfx_\reg$ in \eqref{eq:Tik} by applying rLSMR to the matrices appearing in the leftmost problem in \eqref{eq:Tik} (i.e., following an approach similar to the one described for hybrid rCGLS).

\setlength{\tabcolsep}{6pt} 
\renewcommand{\arraystretch}{1.2} 
\begin{table}
\caption{Summary of the new hybrid randomized methods, as compactly described in \eqref{eq:summary}.}\label{tab:solversummary}
\footnotesize
    \centering
    \begin{tabular}{|c|c|c|c|c|}\hline
         \textbf{method} & \textbf{definition} & $\calR(\bfV_k)$ & $\bfP_k$ & $\bfc_k$\\\hline
hybrid rGMRES  & \eqref{eq: randreg} & $\calK_k(\bfA,\bfb)$ & $\bfH_k\t\bfH_k$ & $\bfH_k\t\beta\bfe_1$\\\hline
hybrid rLSQR  & \eqref{eq:hybrid_rLSQR1}-\eqref{eq:hybrid_rLSQR2} & $\calK_k(\bfA\t\bfA,\bfA\t\bfb)$  & $\bfM_k\t\bfM_k$ & $\bfM_k\t\beta\bfe_1$ \\\hline
hybrid rCGLS & \eqref{eq:hybrid_rCGLS1}-\eqref{eq:hybrid_rCGLS2} & $\calK_k(\bfA\t\bfA,\bfA\t\bfb)$  & $\widetilde{\bfT}_{k+1}\bfM_k$ & $\beta t_{1,1}\bfe_1$ \\\hline
hybrid rLSMR & \eqref{eq:hybrid_rLSMR1}-\eqref{eq:hybrid_rLSMR2} & $\calK_k(\bfA\t\bfA,\bfA\t\bfb)$  & $\!\!\!({\bfT}_{k+1}\bfM_k)\t\!\!({\bfT}_{k+1}\bfM_k)\!\!\!$ & $\!\!\beta t_{1,1}\!(\bfT_{k+1}\bfM_k)\t\!\!\bfe_1\!\!\!\!$ \\\hline
    \end{tabular}
\end{table}
\normalsize

All of the hybrid solvers developed in this section for the Tikhonov-regularized problem compute solutions that belong to approximation subspaces that are independent of the Tikhonov regularization parameter $\lambda$. This property allows us to cheaply approximate different solutions for different values of $\lambda$ and/or to adaptively set $\lambda$ at every iteration by solving the associated projected regularized problems for different values of $\lambda$ without having to recompute the corresponding approximation subspace; as far as the number of iterations $k\leq K\ll \min(m,n)$, the computational cost of doing this is negligible. Many parameter choice rules commonly used to set $\lambda$ for the full-dimensional Tikhonov problem have been successfully applied (with appropriate modifications) to adaptively set $\lambda$ at each iteration of standard (deterministic) hybrid methods; see \cite{chung2024computational, gazzola2020survey, gazzola2015survey} and the references therein. Here we focus on the discrepancy principle and (weighted) GCV only, and briefly describe how they should be modified for use together with the new hybrid randomized solvers. In order to highlight the dependence of the solution $\bfx_k$ on $\lambda$ and to establish a common framework for all the considered hybrid randomized formulations, we use the notation,
\begin{equation}\label{eq:summary}
\bfx_k(\lambda)=\bfV_k\bfz_k(\lambda)=\bfV_k(\bfP_k+\lambda^2\bfI_k)^{-1}\bfc_k,
\end{equation}
where a summary of the new hybrid randomized methods and the corresponding definitions of $\bfV_k$, $\bfP_k$ and $\bfc_k$ are provided in Table \ref{tab:solversummary}.

The discrepancy principle can be applied if a good estimate of $\|\bfe\|$ (the norm of the noise affecting the data in \eqref{eq:inverseproblem}) is available. For hybrid randomized methods, at the $k$th iteration, $k\leq K$, we find the $\lambda$ that satisfies the nonlinear equation,
\begin{equation}\label{eq:dp}
\tau\|\bfe\|=\|\bfThetanp(\bfA\bfx_k(\lambda) - \bfb)\|=\begin{cases}
\|\bfH_k\bfz_k(\lambda)-\beta\bfe_1\| & \text{for hybrid rGMRES},\\
\|\bfM_k\bfz_k(\lambda)-\beta\bfe_1\| & \text{otherwise},
\end{cases}
\end{equation}
by applying a zero finder. The constant $\tau>1$ appearing in the above equation is a safety factor that accounts for possible inaccuracies in the estimation of $\|\bfe\|$ and the fact that the data fit is evaluated in the sketched seminorm.

The weighted generalized cross-validation (WGCV) method can be applied when the noise level is unknown. At the $k$th iteration, $k\leq K$, of the considered hybrid randomized methods, WGCV prescribes to minimize, with respect to $\lambda$, the functional
\begin{equation}\label{eq:GCV}
G_{w}(\lambda)=\frac{n\|\bfThetanp(\bfb-\bfA\bfx_k(\lambda))\|^2}{(\mbox{trace}(\bfI_m-w\bfA\bfA_\reg^\dagger(\lambda)))^2}\,, \quad\mbox{where $w>0$ is a weight.}
\end{equation}
The numerator of $G_w(\lambda)$ can be expressed as in \eqref{eq:dp}; the matrix $\bfA_\reg^\dagger(\lambda)$ appearing in the denominator of $G_w(\lambda)$ is such that 
$\bfx_k(\lambda)=\bfA_\reg^\dagger(\lambda)\bfb$ and is specific of each hybrid method. Namely,
\[
\bfA\bfA_\reg^\dagger(\lambda)=\begin{cases}
\bfV_{k+1}\bfH_k(\bfP_k+\lambda^2\bfI_k)^{-1}\bfH_k\t\bfV_{k+1}\t(\bfThetan)\t\bfThetan,\\
\bfU_{k+1}\bfM_k(\bfP_k+\lambda^2\bfI_k)^{-1}\bfM_k\t\bfU_{k+1}\t(\bfThetanp)\t\bfThetanp,\\
\bfU_{k+1}\bfM_k(\bfP_k+\lambda^2\bfI_k)^{-1}\widetilde{\bfT}_{k+1}\bfU_{k+1}\t(\bfThetanp)\t\bfThetanp,\\
\bfU_{k+1}\bfM_k(\bfP_k+\lambda^2\bfI_k)^{-1} (\bfT_{k+1}{\bfM}_{k})\t  \bfT_{k+1} \bfU_{k+1}\t(\bfThetanp)\t\bfThetanp,\\
\end{cases}
\]
for hybrid rGMRES, rLSQR, rCGLS and rLSMR, respectively. For the denominator of \cref{eq:GCV}, the trace of these matrices can be computed efficiently using the cyclic property.
The weight $w$ in \eqref{eq:GCV} can be adaptively set as described in \cite{ChNaOL08}; choosing $w=1$ in \eqref{eq:GCV} results in the (standard) GCV.

We conclude this section by emphasizing that the functionals in \eqref{eq:dp} and \eqref{eq:GCV} can be efficiently evaluated for different values of $\lambda$, with negligible computational overhead. 


\section{Numerical Results}
\label{sec:numerics}
In this section, we illustrate the effectiveness of the randomized Krylov methods and their hybrid variants for solving inverse problems, compared to their non-randomized counterparts. We judge the reconstruction quality by the $l_2$-norm relative reconstruction error, i.e., by computing $\|\bfx_\rec - \bfx_\true\|/\|\bfx_\true\|$, where $\bfx_\rec$ is the reconstruction computed by each solver. We consider two test problems, generated via \textsc{IR Tools} \cite{2019GazzolaIRtools}: namely, an image deblurring problem and a seismic tomography problem. Moreover, for Tikhonov regularization, we consider randomized methods with fixed $\lambda$ as well as hybrid randomized methods with automatic selection of the regularization parameter.
When naming the methods, all randomized approaches have an `r' in front. All hybrid methods are denoted with a dash and appended parameter or parameter choice method.  For example, rGMRES-$\lambda$ is rGMRES with a fixed regularization parameter at each iteration, and rGMRES-opt is rGMRES where the optimal regularization parameter is selected at each iteration (i.e., the one minimizing the reconstruction error at each iteration). Other combinations include rGMRES-gcv and rGMRES-dp, where the GCV and the discrepancy principle are applied to select the regularization parameter at each iteration, respectively. LSQR and rLSQR can replace GMRES and rGMRES in the acronyms listed so far. Our implementations of the methods described in this paper, together with a selection of the tests illustrated in this section, will be available at \url{https://github.com/juliannechung/randKrylov}.

\subsection{Experiment 1: rGMRES and its hybrid version}
\begin{figure}[bthp]
    \centering
    \includegraphics[width=.8\textwidth]{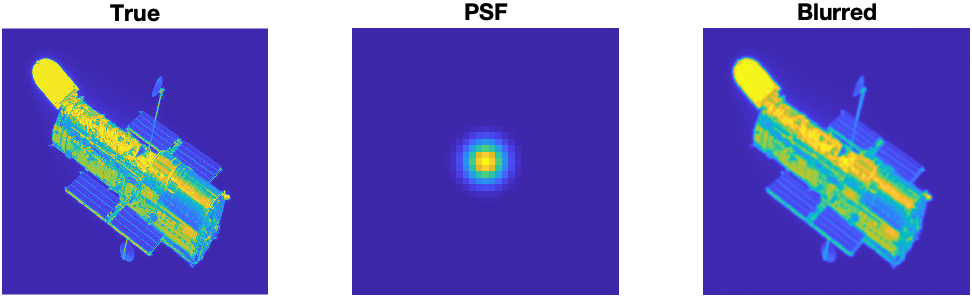}
    \caption{Image debluring example. From left to right: true image, point spread function (PSF), and observed, blurred image.}
    \label{fig:ex_RGMRES}
\end{figure}

\begin{figure}[bthp]
    \centering
    \includegraphics[width=.8\textwidth]{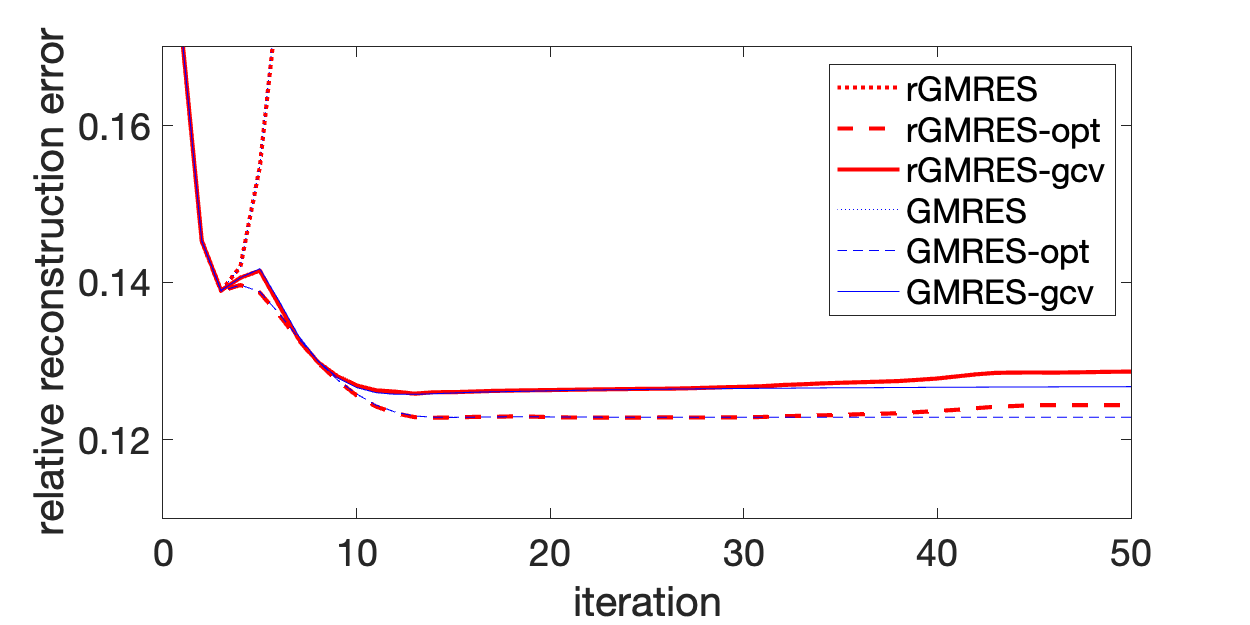}
    \caption{Image deblurring example, with sampling parameter $\ell_n=13108$ (5\% of the size of the matrix $\bfA$). Relative reconstruction error norms per iteration for rGMRES, along with the hybrid rGMRES method with the optimal regularization parameter and the GCV selected regularization parameter at each iteration. Analogous results for standard GMRES and hybrid GMRES are provided for comparison.}
    \label{fig:err_RGMRES}
\end{figure}

Here we consider an example from image deblurring and investigate the performance of rGMRES and hybrid rGMRES. The deblurring example is given in \cref{fig:ex_RGMRES}, where the goal is to reconstruct the true image of size 512 $\times$ 512 pixels (displayed in the leftmost frame), which was corrupted by mild blur defined by the point spread function (PSF, displayed in the middle frame), and additive Gaussian noise, so that the observation (displayed in the rightmost frame) contains $1\%$ noise level. 
Since the forward matrix $\bfA$ is square of order $262144$, we compare the rGMRES and hybrid rGMRES methods to the GMRES and hybrid GMRES methods for this problem. We perform at most $K=50$ iterations of these solvers. 
The relative reconstruction error norms 
per iteration are displayed in Figure \ref{fig:err_RGMRES}. We observe that both GMRES and rGMRES exhibit the so-called semiconvergence behavior, where the relative error norm decreases in early iterations but then increases rapidly.  rGMRES-opt and rGMRES-gcv follow GMRES-opt and GMRES-gcv respectively. Although, in later iterations, the randomized approaches begin to deviate from the nonrandomized counterparts, the error is small.  
This difference is not noticeable from the approximate relative residual norms (not shown here).  
For all of these experiments, we use the SHRT with $\ell_n = 13106.$ This is 5\% of the vector dimension $n$.

\begin{figure}[bthp]
    \centering
    \includegraphics[width=.48\linewidth]{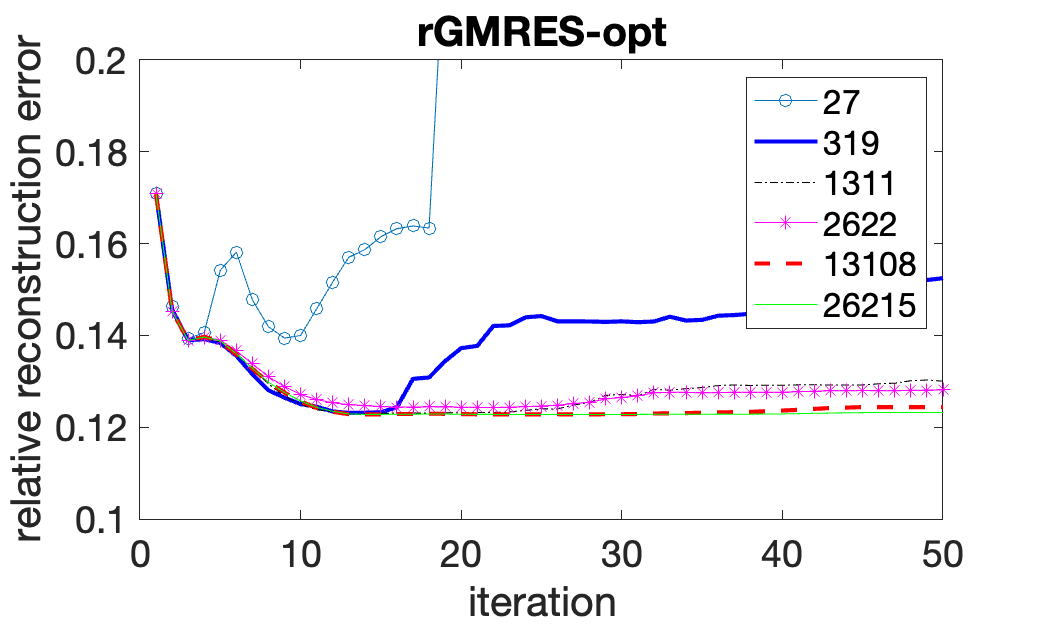}
    \includegraphics[width=.48\linewidth]{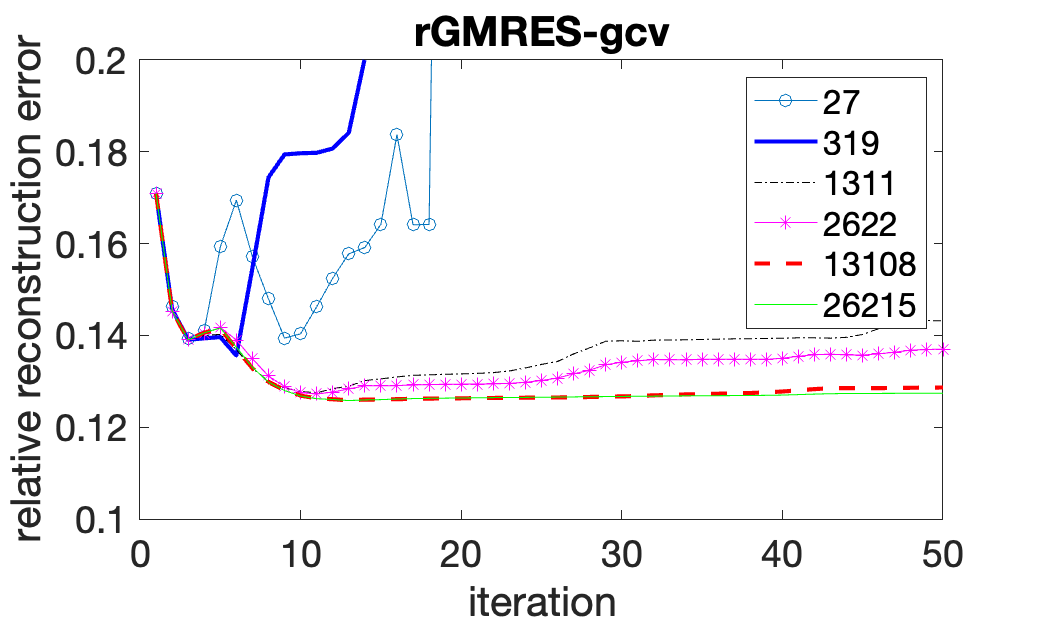}
    \caption{Image deblurring example. For hybrid rGMRES equipped with the optimal regularization parameter and the GCV computed parameter at each iteration, we compare results for different sample sizes $\ell_n$ (noted in the legend).  In particular, the red dashed line corresponds to the result in Figure \ref{fig:err_RGMRES}; the blue solid line corresponds to the default value \eqref{ellnesth}.}
    \label{fig:recon_RGMRES}
\end{figure}

Next, we investigate the impact of $\ell_n$ on the rGMRES results.  Specifically, in \cref{fig:recon_RGMRES}, we consider different values for $\ell_n$ and provide relative reconstruction error norms per iteration for the case where the optimal regularization parameter (left plot) and the GCV computed parameter (right plot) are selected at each iteration.  For the former we observe that, if $\ell_n$ is too small, since we can expect poor inner product approximations during the iterative method, the error behavior is unstable. For large values of $\ell_n$, the algorithm performs similarly to GMRES-opt (blue dashed line in \cref{fig:err_RGMRES}).  Nevertheless, there are a range of values for $\ell_n$ that give reasonable results.  We note that $\ell_n = 319$ corresponds to the default value \eqref{ellnesth}, also used in \cite{BalabanovRGMRES2022}. The other values for $\ell_n$ are selected to be $0.01\%, 0.5\%, 1\%, 5\%$ and $10\%$ of the dimension $n$. Similar results are observed for rGMRES-GCV (shown in the right plot of \cref{fig:err_RGMRES}) and rGMRES without Tikhonov regularization (not shown due to early semiconvergence behavior). 

Although the GMRES and the randomized GMRES approaches only work for square problems, they can be useful for inverse problems where $\bfA\t$ is not accessible or where $\bfA$ is normal and close to being symmetric; see \cite{JeHa07}. In the randomized context, we have shown that rGMRES and hybrid variants can be used to solve large-scale inverse problems, where the main benefit compared to standard methods is that sketched inner products are replaced by approximate inner products. We observe that for very small $\ell_n$ (relative to the size of the problem), reconstruction errors in the early iterations of rGMRES follow the error curve for GMRES, and that for hybrid rGMRES approaches, the regularization parameter can be selected using a standard method (e.g., GCV).  Further numerical investigations for rGK based solvers are provided in the next section.

\subsection{Experiment 2: rGK-based solvers and their hybrid versions}
\begin{figure}[bthp]
    \centering
    \includegraphics[width=.5\textwidth]{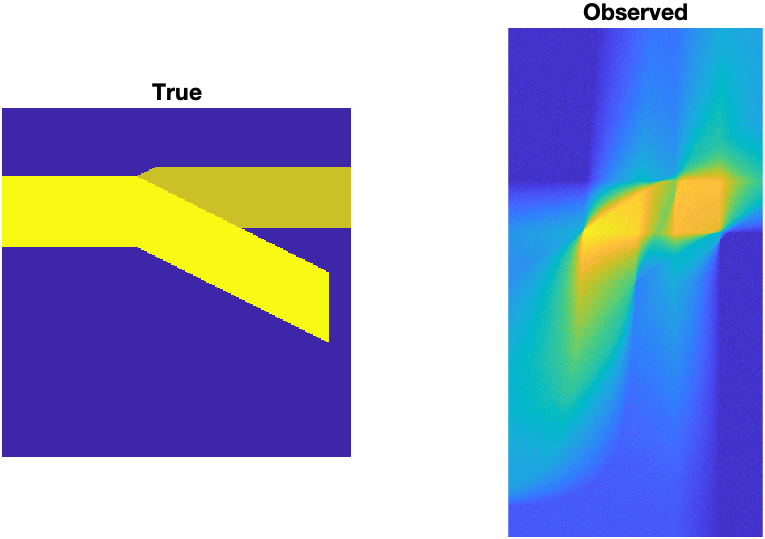}
    \caption{Seismic tomography example. From left to right, true image and observations.}
    \label{fig:seismic}
\end{figure}

The second experiment consists of a seismic travel-time tomography example, where the aim is to reconstruct an image of size 256 $\times$ 256 pixels representing slowness (the reciprocal of the sound speed), from observed data collected along waves modeled as straight lines. Noise was added to the data, in such a way that the observation contained a $4\%$ noise level. The true image and observations are provided in \cref{fig:seismic}. The corresponding matrix $\bfA$ is of size $131,072 \times 65,536$.

\begin{figure}[bthp]
    \centering
    \includegraphics[width=.49\linewidth]{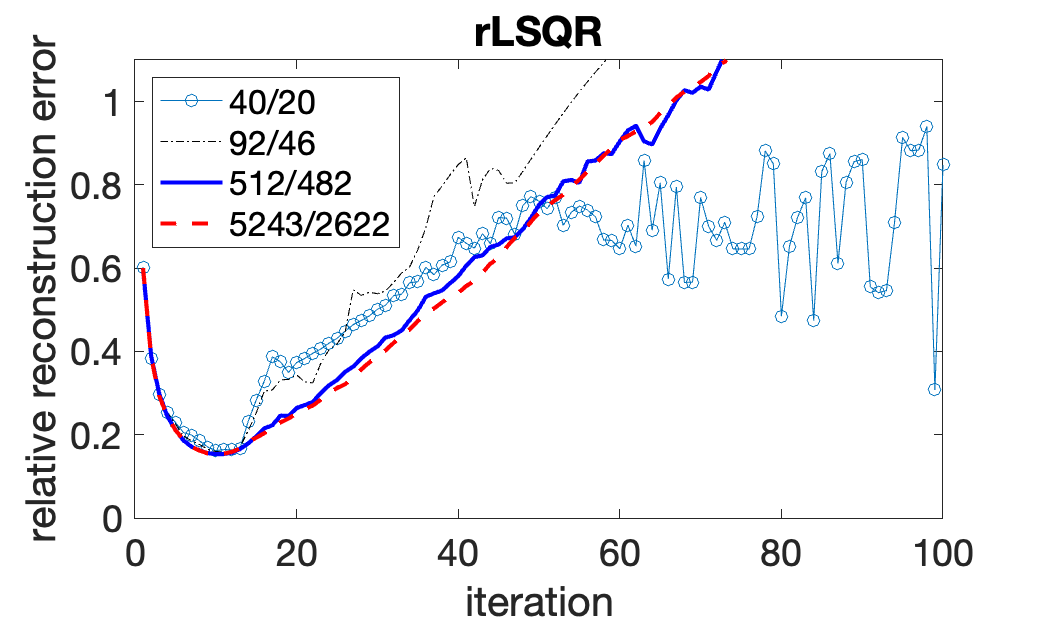}
    \includegraphics[width=.49\linewidth]{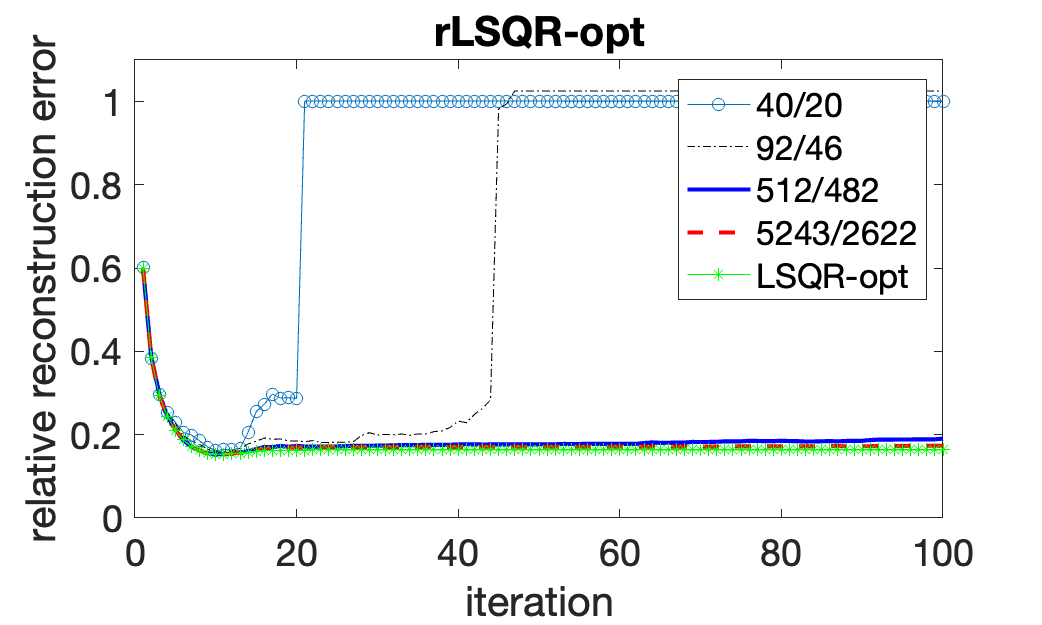}
    \caption{Seismic tomography example. Relative reconstruction error norms per iteration for different sample sizes $\ell_m/\ell_n$ (noted in the legend).  The left plot shows results for rLSQR and the right plot shows results for hybrid rLSQR with the optimal regularization parameter selected at each iteration.}
    \label{fig:seismic_ksamp}
\end{figure}

\begin{figure}[bthp]
    \centering
    \includegraphics[width=.8\textwidth]{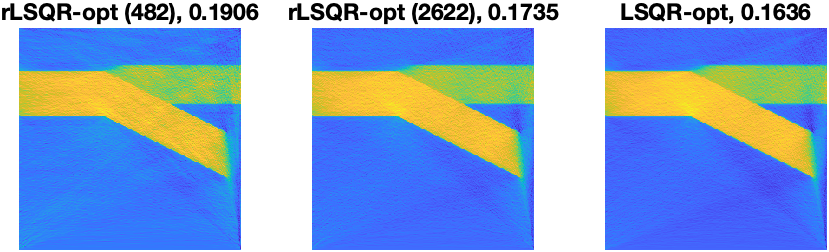}
    \caption{Seismic tomography example. Left and middle frames: reconstructions after $100$ iterations of rLSQR-opt with different sample sizes taken among the ones tested in \cref{fig:seismic_ksamp}; the value of $\ell_n$ is provided in parentheses. Right frame: reconstruction after $100$ iterations of LSQR-opt.  The relative reconstruction error for each image is provided.}
    \label{fig:seismic_recon}
\end{figure}

First we investigate the performance of rLSQR with different sample sizes. The maximum number of iterations is set to $100$, so the sampling parameters automatically selected according to \eqref{ellnesth} are $\ell_m = 512$ and $\ell_n = 482$. We consider sample sizes of $\ell_m/\ell_n = 40/20$, $92/46$, and $5243/2622$, corresponding to $0.03\%$, $0.07\%$ and $4\%$ of $m/n$, respectively. The relative reconstruction error norms per iteration are provided in the left frame of Figure \ref{fig:seismic_ksamp}. Notice the expected semiconvergence behavior in all plots, where early iterations perform similarly, even for very small sample sizes. However, for very small sample sizes, errors increase quickly and exhibit erratic behavior. Somewhat similarly, for the hybrid rLSQR method with the optimal regularization parameter selected at each iteration, in the right frame of Figure \ref{fig:seismic_ksamp} we see that  early reconstructions are good even for small sample sizes, but errors start to increase quickly after, and eventually, minimal improvements are made. 
On the other hand, for a sufficiently large sample size, reconstruction errors remain small for all iterations: this behavior mimics the hybrid LSQR method with optimal regularization parameter (LSQR-opt), which is provided for comparison. For a visual comparison, in Figure \ref{fig:seismic_recon} we provide the reconstruction from rLSQR-opt at 100 iterations for $\ell_n =482$ and $\ell_n=2622$, along with the LSQR-opt reconstruction.  We observe that, as anticipated, some details are lost in the rLSQR-opt reconstructions, especially for small sample parameters, but both reconstructions are still reasonable.

Next, we investigate the effect of the choice of the regularization parameter $\lambda$ on the hybrid rGK-based solvers. For these experiments, we allow the automatically selected sample size $\ell_m/ \ell_n = 512/482$. 

\begin{figure}[bthp]
    \centering
    \includegraphics[width=.9\linewidth]{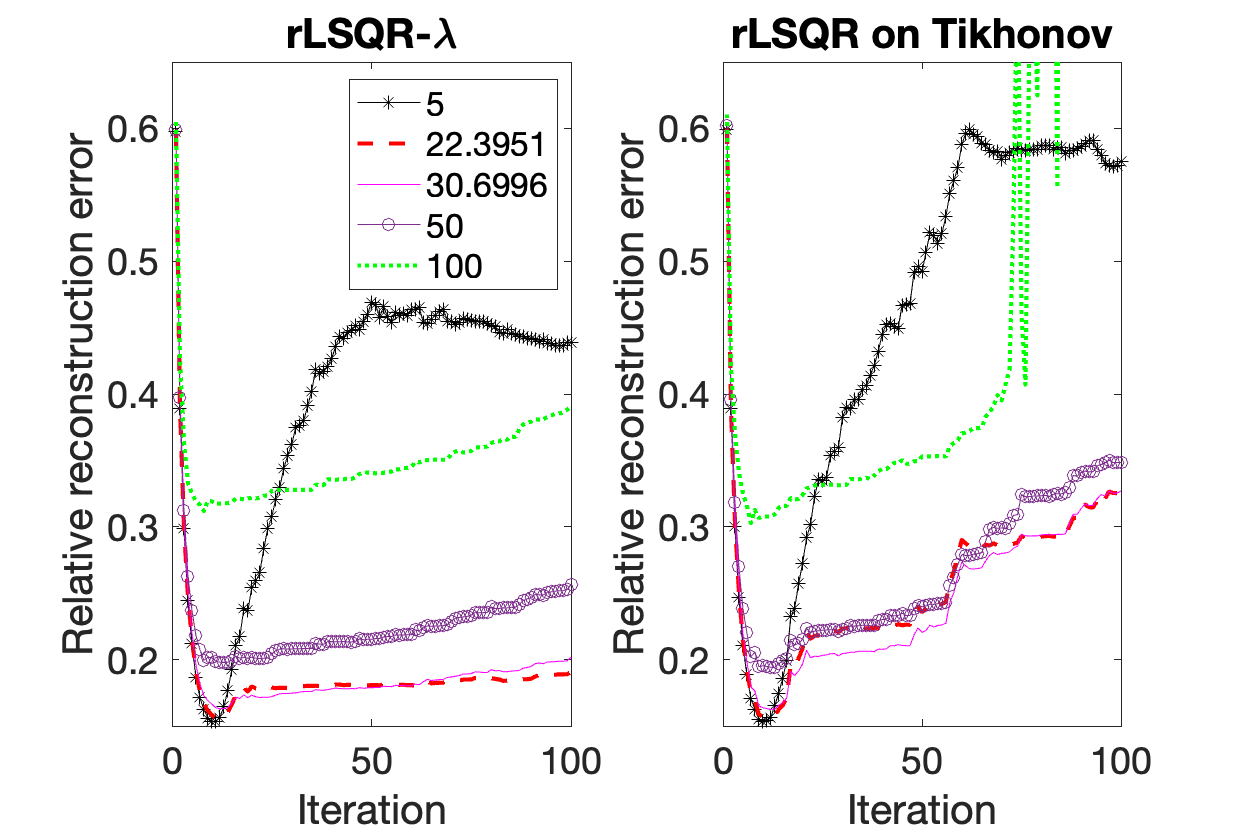}
    \caption{Seismic tomography example. Left frame: comparison of relative reconstruction errors for hybrid rLSQR \eqref{eq:hybrid_rLSQR1}-\eqref{eq:hybrid_rLSQR2}, with a fixed regularization parameter (its value is noted in the legend).  Right frame: relative reconstruction errors for rLSQR applied to the Tikhonov problem as in \eqref{eq:rLSQR_dampTik}, where the regularization parameter is fixed in advance. The legend is the same for both plots.}
    \label{fig:seismic_regparam}
\end{figure}

\begin{figure}[bthp]
    \centering
    \includegraphics[width=.7\linewidth]{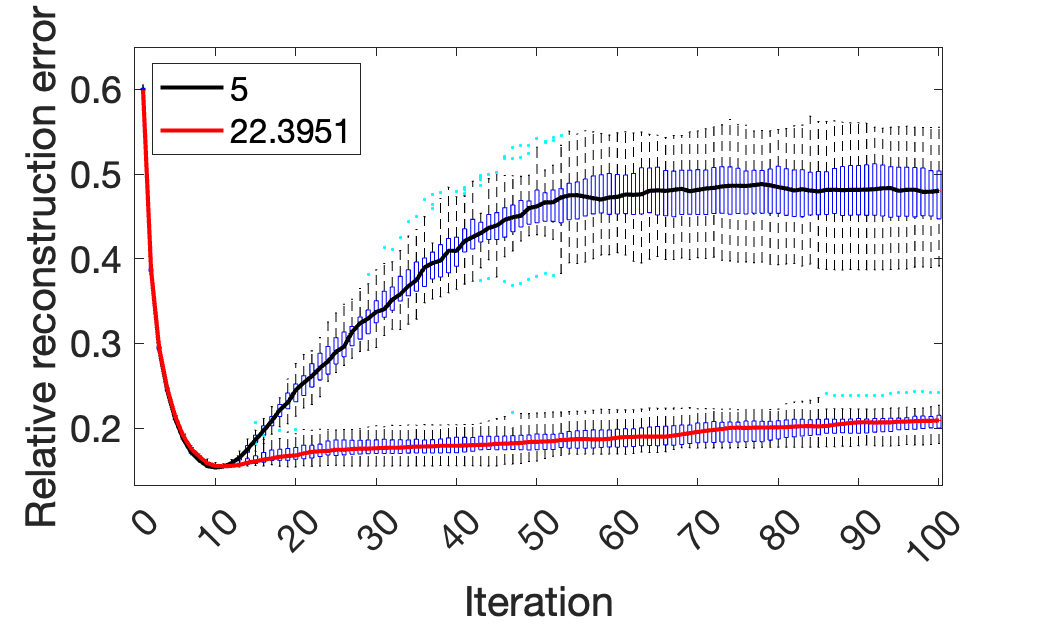}
    \caption{Seismic tomography example. Statistics for the rLSQR-$\lambda$ relative reconstruction errors, with fixed regularization parameter, over $24$ runs. The solid line corresponds to the median values and the edges of the boxes indicate the $25$th and $75$th percentiles. Outliers are plotted separately. }
    \label{fig:seismic_box}
\end{figure}

We first consider the case where a fixed regularization parameter is used at each iteration. In the left frame of Figure \ref{fig:seismic_regparam} we provide relative reconstruction errors for rLSQR-$\lambda$ for different values of $\lambda$.  We note that $22.3951$ is the optimal regularization parameter for rLSQR-opt at iteration $100$, and $30.6996$ is the optimal regularization parameter for LSQR-opt at iteration $100$. For regularization parameters that are too small, we observe semiconvergence, and for regularization parameters that are too large we obtain overly regularized solutions (again corresponding to large errors). Still in Figure \ref{fig:seismic_regparam}, in the right frame, for comparison, we investigate how rLSQR performs when applied to the Tikhonov problem, for a fixed regularization parameter and as described in \cref{eq:rLSQR_dampTik}. We use the same set of regularization parameters as in the left plot of \Cref{fig:seismic_regparam}. Note that, since the augmented matrix in \cref{eq:rLSQR_dampTik} has $m+n$ rows, the sample sizes selected according to \eqref{ellnesth} are $\ell_{m+n}/\ell_n = 530/482$. We observe that applying rLSQR to the Tikhonov problem results in more irregular behavior and overall larger errors than rLSQR-$\lambda$. 
%
To back up the results displayed in the left frame of \cref{fig:seismic_regparam}, since there is randomness involved in the sketching process, we run the same solver $24$ times for two fixed values of the regularization parameter (among the ones considered in \Cref{fig:seismic_regparam}) and, in Figure \ref{fig:seismic_box}, we provide box plots for the relative reconstruction errors at each iteration, highlighting median and the $25$th and $75$th percentiles. Outliers are plotted separately. We see that randomness in the sketching can affect the results, but the general trends remain the same.

\begin{figure}[bthp]
    \centering
    \includegraphics[width=\linewidth]{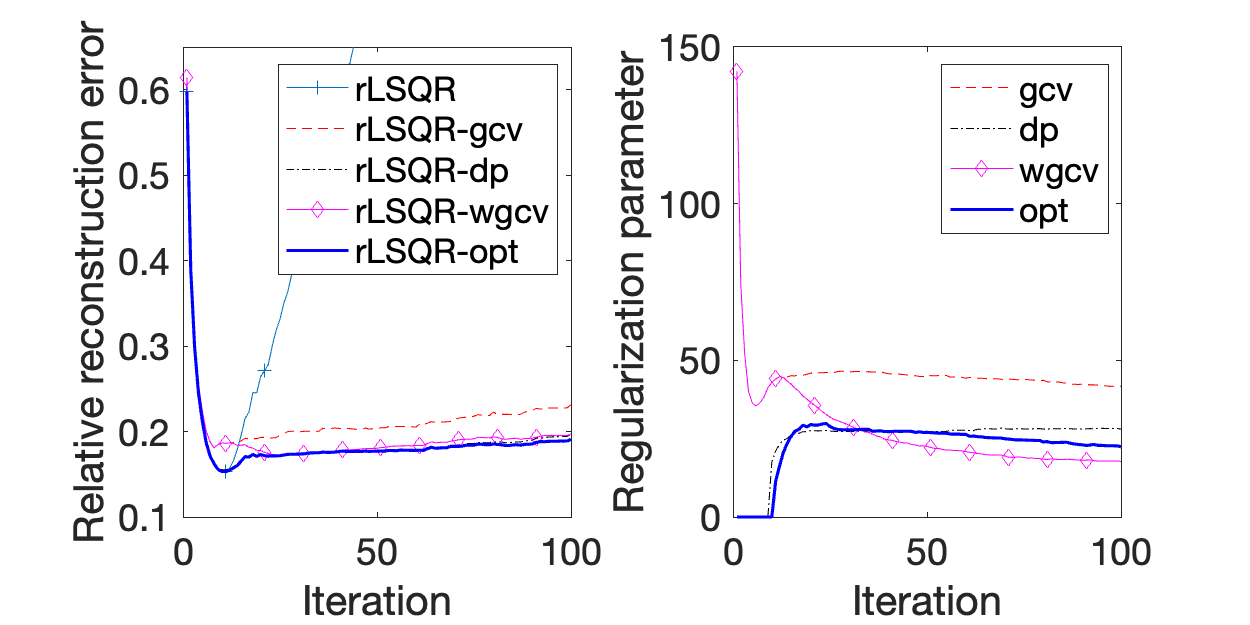}
    \caption{Seismic tomography example. Left frame: comparison of relative reconstruction errors per iteration for hybrid rLSQR with different regularization parameter selection methods. Regularization parameters selected at each iteration are provided in the right plot.}
    \label{fig:err_RLSQR}
\end{figure}

Next, we consider different adaptive regularization parameter selection methods. In the left frame of \cref{fig:err_RLSQR} we show the relative reconstruction errors per iteration of the hybrid rLSQR method. Namely, we provide results for GCV, DP, WGCV, and when the optimal regularization parameter is selected at each iteration. The corresponding regularization parameters computed at each iteration are provided in the right frame of \cref{fig:err_RLSQR}.  Note that the optimal regularization parameter is zero in early iterations, since the problem is still well-conditioned and no regularization is needed. The discrepancy principle and the WGCV method provide good approximations of the regularization parameter at later iterations (i.e., close to its optimal values), while the GCV method overestimates the regularization parameter.

\begin{figure}
    \centering
     \includegraphics[width=\linewidth]{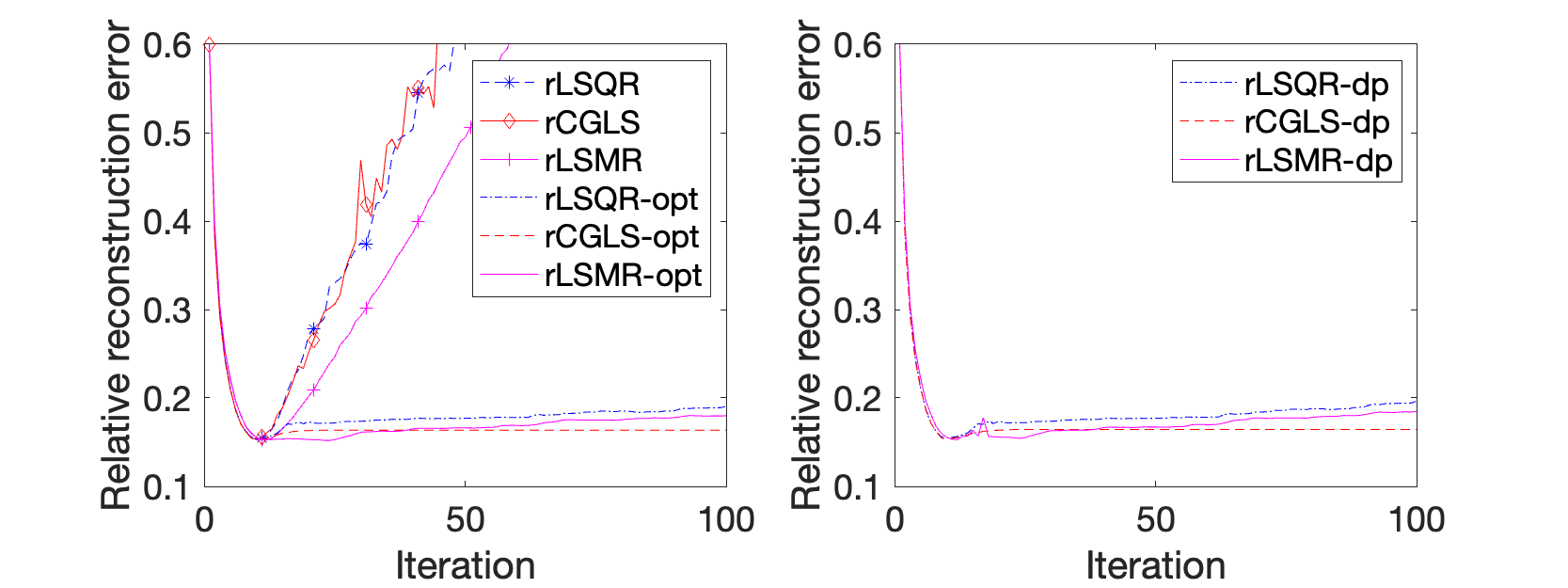}
    \caption{Seismic tomography example. Comparison of different rGK-based solvers (i.e., rLSQR, rCGLS, and rLSMR) for the seismic example.  The left plot contains relative reconstruction error norms per iteration for rGK-based methods with no regularization term and the hybrid variant with the optimal regularization parameter selected at each iteration.  In the right plot we provide the relative reconstruction errors for the hybrid variant with the regularization parameter selected by the discrepancy principle.}
    \label{fig:rGKsolvers}
\end{figure}

Finally, we investigate the performance of the different solvers based on the rGK algorithm. In the left plot of \cref{fig:rGKsolvers} we provide relative reconstruction errors per iteration for rLSQR, rCGLS, rLSMR and their hybrid versions with the optimal regularization parameter selected at each iteration. In the right frame of \cref{fig:rGKsolvers}, we use the regularization parameter selected by the discrepancy principle at each iteration. We observe that all of the considered methods exhibit similar behavior, with the rLSMR methods having slightly smaller reconstruction errors. 

\begin{figure}[bthp]
    \centering
    \includegraphics[width=\linewidth]{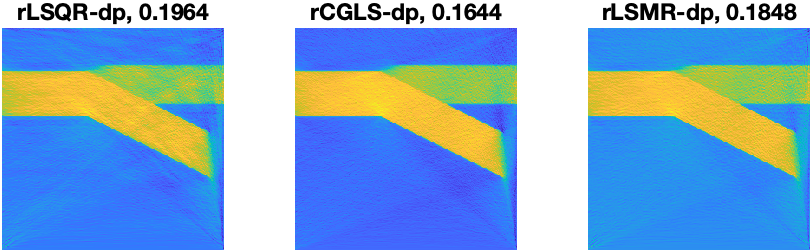}
    \caption{Seismic tomography example. Reconstructions after $100$ iterations of rLSQR-dp, rCGLS-dp, and rLSMR-dp.}
    \label{fig:seismic_recon_rGK}
\end{figure}

With the seismic tomography example, we have shown that rGK-based solvers can be used to solve inverse problems with rectangular discretized forward operator $\bfA$, where regularized solutions can be obtained with early termination (e.g., before the errors increase) or via a hybrid approach, where various parameter selection methods can be included. Compared to standard Krylov methods, all of the randomized Krylov approaches (including the ones tested for the image deblurring example) require the selection of a sampling parameter. For problems where only a few iterations are needed to obtain a reasonable solution (e.g., problems with very high noise levels or severe ill-conditioning), smaller values of the sampling parameter may be sufficient. Choosing a good regularization parameter is also important, although it may not be able to compensate for errors introduced in the basis generation.




\section{Conclusions}
\label{sec:conclusions}
In this work, we develop randomized Krylov methods for large-scale inverse problems, where sketched inner products are used within a randomized Gram-Schmidt process to generate suitable projected problems. For square problems, we consider randomized GMRES and introduce its hybrid variant for solving Tikhonov-regularized problems. For rectangular matrices, we develop a randomized Golub-Kahan approach, where each iteration requires a matrix-vector multiplication with $\bfA$ and $\bfA\t$, and all the inner products are approximated with sketched inner products.  We describe new iterative solvers based on the the randomized Golub-Kahan method, and describe hybrid variants of all approaches, with automatic Tikhonov regularization parameter selection at each iteration. Numerical results show that randomized Krylov methods can be used to solve inverse problems, taking into account that semiconvergence appears when using randomized GMRES and LSQR, and when the sample size for the randomized inner products approximation is too small. In the latter case, even a good regularization parameter may not be able to provide a good reconstruction if not terminated early. 

There are various promising directions for future work. For example, there are many choices for selecting the sketching matrix, and exploiting any matrix structure may improve the results and significantly improve run time. We are currently working on a highly efficient computational implementation and an analysis of the computational complexity of the algorithm.  Moreover, there are potential extensions to settings such as Tikhonov regularization in general form and flexible Krylov methods to handle $p$-norm regularization, $0<p\leq 1$. Finally, further research is needed to understand the theoretical properties of the low-rank matrix approximations associated with the new rGK algorithm, which may eventually be useful for estimating matrix functions or for uncertainty quantification.

\section{Acknowledgments} We would like to thank Leonardo Robol for many helpful and stimulating discussions about randomized numerical linear algebra.

\bibliographystyle{siamplain}
\bibliography{references}

\end{document}